\documentclass{svjour3}
\smartqed  
\usepackage{graphicx}
\usepackage{latexsym,bm,amsmath,amssymb}
\usepackage{rotating}
\usepackage{multirow,booktabs,cases}
\usepackage[misc]{ifsym}
\usepackage[style=1]{mdframed}
\usepackage{algorithm,algorithmic}
\usepackage[colorlinks=true]{hyperref}
\hypersetup{urlcolor=blue, citecolor=blue,linkcolor=blue}
\usepackage{epstopdf}
\usepackage{bm}
\usepackage{threeparttable}
\usepackage{epstopdf}
\usepackage{caption}
\usepackage{subcaption} %
\def\x{{\bf x}}
\def\y{{\bf y}}
\def\z{{\bf z}}
\def\w{{\bf w}}
\def\u{{\bf u}}
\def\v{{\bf v}}
\def\t{{\bf t}}
\def\0{{\bf 0}}

\def\F{{\mathcal F}}
\def\T{{\mathcal T}}
\allowdisplaybreaks

\begin{document}
\graphicspath{{./PIC/}}

\title{Tensor Based Proximal Alternating Minimization Method for A Kind of Inhomogeneous Quartic Optimization Problem}
\titlerunning{A PAM Method for Fourth Degree POP}     

\author{Haibin Chen \and {Yixuan Chen} \and Chunyan Wang \and Qi Fan}


\institute{Authors are listed alphabetically. \\
H. Chen {\Letter} \at School of Management Science, Qufu Normal University, Rizhao, Shandong, China, 276800. \\
\email{chenhaibin508@qfnu.edu.cn}
\and Yixuan. Chen \at School of Management Science, Qufu Normal University, Rizhao, Shandong, China, 276800. \\
\email{yixuanchen1201@163.com}
\and C. Wang \at School of Management Science, Qufu Normal University, Rizhao, Shandong, China, 276800.\\
\email{wchunyanwang562@163.com}
\and Q. Fan \at School of Management Science, Qufu Normal University, Rizhao, Shandong, China, 276800. \\
\email{fanqi1526@163.com}}

\date{Received: date / Accepted: date}

\maketitle

\begin{abstract}
In this paper, we propose an efficient numerical approach for solving a specific type of quartic inhomogeneous polynomial optimization problem inspired by practical applications. The primary contribution of this work lies in establishing an inherent equivalence between the quartic inhomogeneous polynomial optimization problem and a multilinear optimization problem (MOP). This result extends the equivalence between fourth-order homogeneous polynomial optimization and multilinear optimization in the existing literature to the equivalence between fourth-order inhomogeneous polynomial optimization and multilinear optimization. By leveraging the multi-block structure embedded within the MOP, a tensor-based proximal alternating minimization algorithm is proposed to approximate the optimal value of the quartic problem. Under mild assumptions, the convergence of the algorithm is rigorously proven. Finally, the effectiveness of the proposed algorithm is demonstrated through preliminary computational results obtained using synthetic datasets.
\end{abstract}

\section{Introduction}\label{Int}

In this paper, we consider the following inhomogeneous polynomial optimization problem (IPOP) over a spherical constraint:
\begin{equation}\label{e1}
\begin{aligned}
\min_{\textbf{\text{x}} \in \mathbb{R}^n} &~f(\x)=\frac{\theta}{2} \sum_{i=1}^n x_i^4 + \textbf{\text{x}}^\top B \textbf{\text{x}}\\
\textup{s.t.}   &~\|\textbf{\text{x}}\| = 1,
\end{aligned}
\end{equation}
where $\theta > 0$ is a fixed constant and $B$ is a symmetric matrix with positive diagonal entries and nonpositive off-diagonal entries. Note that problem (\ref{e1}) maybe nonconvex.

The main motivation for studying (\ref{e1}) arises from the discretization of energy functional minimization tasks for calculating the ground state of Bose-Einstein Condensate (BEC)  or other electronic structures \cite{2022-7,2022-10}. According to the requirements of different applications, various types of BECs have been developed, including rotating BEC with quantum vortices \cite{BEC-11}, BEC with dipole-dipole interactions \cite{BEC-12}, and multi-component BEC \cite{BEC-13}.
It is generally NP-hard to solve such problems \cite{2022-17,2022-34}, thus we focus on a relatively simple but significant class of applications: identifying the ground state of the non-rotating BEC problem (\ref{e1}). Typically, this state is characterized as the minimizer within the energy functional minimization problem. Through appropriate discretization techniques (refer to Section 4),  the matrix $B$ is formulated to represent  the sum of the discretized negative Laplace operator and a positive diagonal matrix. As we will see later,  even for such a relatively simple BEC problem, the calculation of its numerical ground state solution will involve a quartic nonconvex optimization problem with spherical constraints. For an in-depth understanding and contextual background, the reader is directed to  \cite{2022-3} and the references cited therein.

To the best of our knowledge, the methods for finding the ground state solution of BEC can be divided into two categories.
The first category involves optimization methods such as semidefinite programming method \cite{2022-26,2022-17}, projection gradient method \cite{BEC-22,BEC-35}, preconditioned conjugate gradient method \cite{BEC-36,BEC-37},
Riemann manifold method \cite{BEC-38,2022-34}, regularized Newton method \cite{2022-35} and gradient flow method \cite{2022-3,BEC-30,BEC-26,BEC-27}. The second category includes methods from numerical linear algebra \cite{BEC-31,2022,BEC-4}, where the origin problem is transformed to a series of linear eigenvalue problems.

On the other hand, by the spherical constraint of (\ref{e1}), certain specialized algorithms were provided to solve it efficiently. For instance, Tang {\it et al}. \cite{Chen2021-11,BEC-48} employed the Shifted Symmetric High-Order Power Method to solve (\ref{e1}). By exploiting the multi-block structure of the polynomial functions, Chen {\it et al}. \cite{Chen2021-2} introduced a maximum block improvement (MBI) method, which draws inspiration from the well-known block coordinate descent (BCD) method (e.g., see \cite{Chen2021-24,Chen2021-27}). Considering the equivalent multilinear function in \cite{Chen2021-2}, only one block in MBI method is updated at each iteration. To sufficiently apply the multilinear structure, Jiang {\it et al}. \cite{Chen2021-10} redefined the problem into a constrained linear polynomial optimization model.
By the linear model, it is particularly advantageous for employing the advanced ADMM technique. Instead of the only one block update rule, ADMM updates all blocks in a sequential (or alternating) order with the help of augmented Lagrangian function \cite{Chen2021-10}. To weaken the theoretical requirements of the ADMM in \cite{Chen2021-10}, Wang {\it et al}. proposed a block improvement method (BIM) \cite{Chen2021-29}, and two blocks are updated in a sequential order.
Recently, Huang {\it et al}. \cite{2022-17} studied the nonconvex model (\ref{e1}). The convergence to the global optimum of the inexact ADMM was obtained, and numerical experiments for applications in non-rotating BECs were given to validate the theories \cite{2022-17}.
Above all, we can see that the block updating strategy is highly effective for this kind of problem.

In this paper, we will solve (\ref{e1}) by the sparse tensor composed from the coefficients of the polynomial, and a tensor based proximal alternating minimization (PAM) algorithm will be given with a cyclic rule. To do that, we firstly give a sufficient condition to ensure the equivalence between the IPOP(\ref{e1}) and a mulitlinear optimization problem (MOP). Then, with the help of the MOP, we propose the PAM algorithm, which guarantees all subproblems with analytic optimal solutions. Finally, several numerical examples are given to verify the PAM method.
	
The rest of the paper is organized as follows. In Section 2, we recall some notations and preliminaries on tensors and polynomials. In Section 3,
we first transform (\ref{e1}) into an MOP. Then, we consider an augmented problem over a unit sphere, which shares the same optimal solution with the MOP. Furthermore, the PAM algorithm is given for the augmented MOP. We will also show that the sequence generated by BIM embedded with PAM algorithm converges to a KKT point of the IPOP. In Section 4, several numerical examples are given to support the efficiency of the PAM. Conclusions and some final remarks are given in Section 5.

\section{Prelimimary}\label{Sec2}
In this section, we recall some useful notations and preliminaries.
Let $\mathbb{R}^n$ be the $n$ dimensional real Euclidean space. Generally, we use lowercase letters $a, b$ and  bold lowercase letters ${\bf x}, {\bf y}$ to denote scalars and vectors, respectively. Particularly, let $ \mathbf{0}$ denote the zero vector and let $\mathbf{1}$ denote the all one vector.
Matrices are denoted by capital letters such as $ A, B, C$ and tensors are denoted by calligraphic letters such as $\mathcal{A}, \mathcal{B}, \mathcal{C}$. The superscript $^\top$ stands for the transpose of vectors and matrices.
Assume that $m$, $n$ are integers, and $m, n \geq 2$, unless otherwise stated, and denote the set of all $m$-th order $n$-dimensional real tensors as $\mathbb{R}^{m,n }$. For $\mathcal{A}=(a_{i_1 \dots i_m})\in\mathbb{R}^{m,n} $, if its entries are invariant under any permutation of its indices, then $\mathcal{A}$ is called a symmetric tensor. Denote the set of all $m$-th order $n$-dimensional real symmetric tensors by $\mathbf{S}^{m,n}$.	
Throughout this paper, the identity tensor $\mathcal{I}\in \mathbf{S}^{m,n}$ is given with entries such that $\mathcal{I}_{i_1\cdots i_m}=1$ if $i_1=\cdots=i_m$ and $\mathcal{I}_{i_1\cdots i_m}=0$ otherwise.

The standard Euclidean norm of a vector $\x$ is denoted by $\|\x\|=\sqrt{\langle \x,\x\rangle}$.
On the other hand, the Frobenius norm of tensor $\mathcal{A}=(a_{i_1\ldots i_m}) \in \mathbb{R}^{m,n}$ is the usual Euclidean norm defined as
$$
\|\mathcal{A} \| := \sqrt{\sum_{ i_1, i_2, \cdots,i_m \in [n] } {a_{i_1 i_2 \cdots i_m}}^2},
$$
where $[n]:=\left\{1,2,\cdots, n\right\}$.
	
It should be noted that symmetric tensors and homogeneous polynomials have
a one-to-one relationship. Let $\mathbb{R}_m[\mathbf{x
}]$ denote the set of all real polynomials with
degree at most $m$. Suppose $\mathcal{A} = (a_{i_1 i_2 \ldots i_m}) \in \mathbf{S}^{m,n}$, its
corresponding homogeneous polynomial is denoted by $f(\mathbf{x}) \in \mathbb{R}_m[\mathbf{x}]$, $\mathbf{x} \in \mathbb{R}^n$ such
that
\[
f(\mathbf{x}) = \mathcal{A}\mathbf{x}^m = \langle \mathcal{A}, \underbrace{\mathbf{x} \circ \mathbf{x} \circ \ldots \circ \mathbf{x}
}_{m}\rangle =
\sum_{i_1,\ldots,i_m \in [n]} a_{i_1 i_2 \ldots i_m} x_{i_1} x_{i_2} \ldots x_{i_m},
\]
\noindent where “$\circ$” denotes the outer product. If $\mathcal{A}$ is a fourth-order diagonal tensor with all diagonal entries being one, then
$\mathcal{A}\x^4 = \sum_{i=1}^n x_i^4$.

\section{The inhomogeneous quartic optimization problem with spherical constraint}\label{sec_unit}
Let us start by introducing the following multilinear function
$$
F(\mathbf{x}, \mathbf{y}, \mathbf{z}, \mathbf{w}) = \langle\mathcal{F},\mathbf{x}\circ \mathbf{y}\circ \mathbf{z}\circ \mathbf{w}\rangle = \sum_{1 \leq i, j, s, l \leq n } f_{i j s l} x_i y_j z_s w_l , $$
where $\mathbf{x}, \mathbf{y}, \mathbf{z}, \mathbf{w} \in \mathbb{R}^n$, $ \mathcal{F} = (f_{ijsl} ) \in \mathbf{T}^{4,n}$ is the tensor composed by the coefficients of the polynomial $F(\mathbf{x}, \mathbf{y}, \mathbf{z}, \mathbf{w})$ being its associated multilinear function.
On the other hand, when $\x=\mathbf{y}=\mathbf{z}=\mathbf{w}$, the corresponding homogeneous polynomial is $f(\mathbf{x},\mathbf{x},\mathbf{x},\mathbf{x})=\F\x^4$.

\subsection{Homogenization}

We first introduce two basic definitions, which are useful in the $homogenization$.

\begin{definition}\label{def1}
 \textup{\cite{pena2015}} For any $\mathbf{x}\in \mathbb{R}^n$, denote the mapping $M_d : \mathbb{R}^n \rightarrow \mathbf{S}^{d,n}$ such that
$$
M_d(\mathbf{x}) = \underbrace{\mathbf{x} \circ \mathbf{x} \circ \ldots \circ \mathbf{x}  }_{d}.
 $$
	\end{definition}
		
	\begin{definition}\label{def2}
Let $\mathbb{R}_d[\mathbf{x}] = \left\{ p(\mathbf{x}) \mid \deg p \leq d \right\}$ denote the set of polynomials with dimension $n$ and degree at most $d$. For a general polynomial
\[ f(\mathbf{x}) = \sum_{\alpha \in \mathbb{Z}^n_+} f_{\alpha} \mathbf{x}^{\alpha} \in \mathbb{R}_d[\mathbf{x}], \]
the corresponding tensor $\mathcal{T}_f \in \mathbf{S}^{d,n+1}$ is defined as follows:
\begin{equation}\label{e2}
 (\mathcal{T}_f)_{i_1 i_2 \ldots i_d} = \frac{(d - |\alpha|)! \alpha_1! \ldots \alpha_n!}{d!} f_{\alpha},
\end{equation}
where $\alpha=(\alpha_1,\alpha_2,\cdots,\alpha_n)$ is the unique exponent such that
\[ \mathbf{x}^{\alpha} = x_1^{\alpha_1} x_2^{\alpha_2} \ldots x_n^{\alpha_n} = x_{i_1} x_{i_2} \ldots x_{i_d}, \]
and $|\alpha| = \alpha_1 + \alpha_2 + \ldots + \alpha_n\leq d.$ \end{definition}	

By Definition \ref{def1}, we use $M_d(\tilde{\mathbf{x}}) = M_d((1, \mathbf{x}^\top)^\top) \in \mathbf{S}^{d,n+1}$. Then for any $f(\mathbf{x}) \in \mathbb{R}_d[\mathbf{x}]$, it follows that
$$
f(\mathbf{x}) = \langle \mathcal{T}_f, M_d(\tilde{\mathbf{x}})\rangle.
$$
For the sake of clarity, we give an example to explain the detail of the equation above.
\begin{example}
Suppose that $g(\textbf{\textup{x}}) \in \mathbb{R}_3[\mathbf{x}]$ is a polynomial function such that
\[ g(\textbf{\textup{x}}) = x^3_1 + x^3_2 + x^3_3 + x_1x_3 + x_2x_3 + x_1 + x_2 - 2. \]
Let $\mathcal{T} := \mathcal{T}_g \in \mathbf{S}^{3,4}$ and $\mathcal{M} := M_d(\tilde{\mathbf{x}}) \in \mathbf{S}^{3,4}$. Then for any $\textbf{\textup{x}} \in \mathbb{R}^3$, it holds that $g(\mathbf{x}) =
\left\langle \mathcal{T}, \mathcal{M} \right\rangle$, where
	\begin{align*}
&\mathcal{T}_{0,:,:} =
\begin{pmatrix}
-2 &~\frac{1}{3} &~ \frac{1}{3} &~ 0 \\
\frac{1}{3} &~ 0 &~ 0 &~ \frac{1}{6} \\
\frac{1}{3} &~ 0 &~ 0 &~ \frac{1}{6} \\
0 &~\frac{1}{6} &~ \frac{1}{6} &~ 0
\end{pmatrix},
\quad
&\mathcal{T}_{1,:,:} =
\begin{pmatrix}
\frac{1}{3} &~ 0 &~ 0 &~ \frac{1}{6} \\
0 &~ 1 &~ 0 &~ 0 \\
0 &~ 0 &~ 0 &~ 0 \\
\frac{1}{6} &~ 0 &~ 0 &~ 0
\end{pmatrix},\\
&\mathcal{T}_{2,:,:} =
\begin{pmatrix}
\frac{1}{3} &0 & 0 & \frac{1}{6} \\
0 & 0 & 0 & 0 \\
0 & 0 & 1 & 0 \\
\frac{1}{6} & 0 & 0 &0
\end{pmatrix},
\quad
&\mathcal{T}_{3,:,:} =
\begin{pmatrix}
0 & \frac{1}{6} & \frac{1}{6} & 0 \\
\frac{1}{6} & 0 & 0 & 0 \\
\frac{1}{6} & 0 & 0 & 0 \\
0 & 0 & 0 & 1
\end{pmatrix},
	\end{align*}
and
	\begin{align*}
	&\mathcal{M}_{0,:,:} =
	\begin{pmatrix}
	1 & x_1 & x_2 & x_3 \\
	x_1 & x_1^2  & x_1x_2 & x_1x_3 \\
	x_2 & x_2x_1 & x_2^2 & x_2x_3 \\
	x_3 & x_3x_1 & x_3x_2 & x_3^2\\
	\end{pmatrix},
	\quad
	&\mathcal{M}_{1,:,:} =
	\begin{pmatrix}
	x_1 & x_1^2 & x_1x_2 & x_1x_3 \\
    x_1^2 & x_1^3  & x_1^2x_2 & x_1^2x_3 \\
    x_1x_2 & x_1^2x_2 & x_1x_2^2 & x_1x_2x_3 \\
    x_1x_3 & x_1^2x_3 & x_1x_2x_3 & x_1x_3^2\\
	\end{pmatrix},\\
	&\mathcal{M}_{2,:,:} =
	\begin{pmatrix}
	x_2 & x_2x_1 & x_2^2 & x_2x_3 \\
    x_2x_1 & x_2x_1^2  & x_2^2x_1 & x_2x_1x_3 \\
    x_2^2 & x_2^2x_1 & x_2^3 & x_2^2x_3 \\
    x_2x_3 & x_2x_3x_1 & x_2^2x_3 & x_2x_3^2\\
	\end{pmatrix},
	\quad
	&\mathcal{M}_{3,:,:} =
	\begin{pmatrix}
	x_3 & x_3x_1 & x_3x_2 & x_3^2 \\
    x_3x_1 & x_3x_1^2  & x_3x_1x_2 & x_3^2x_1 \\
    x_3x_2 & x_3x_2x_1 & x_3x_2^2 & x_3^2x_2 \\
    x_3^2 & x_3^2x_1 & x_3^2x_2 & x_3^3\\
	\end{pmatrix}.
	\end{align*}
\end{example}

By the above symbols, for the problem (\ref{e1}) with matrix:
$$B=
\begin{pmatrix}
b_{11} & b_{12} & \cdots & b_{1n} \\
b_{12} & b_{22} & \cdots & b_{2n} \\
\vdots & \vdots & \ddots & \vdots \\
b_{1n} & b_{2n} & \cdots & b_{nn}
\end{pmatrix},
$$
then we have
\begin{align*}
f(\mathbf{x})=&\frac{\theta}{2} x_1^4 +\frac{\theta}{2} x_2^4 +\cdots +\frac{\theta}{2} x_n^4 + b_{11} x_1^2+2b_{12}x_1x_2+\cdots+2b_{1n}x_1x_n\\
&+ b_{22}  x_2^2+\cdots +2b_{2n}x_2x_n+\cdots+b_{nn} x_n^2\\
&=\langle \mathcal{T}_f, M_d(\tilde{\mathbf{x}})\rangle,
\end{align*}
where $\tilde{\mathbf{x}}=(1, \mathbf{x}^\top)^\top $, and $\mathcal{T}_f \in \mathbf{S}^{4,n+1}$ is the tensor composed by the coefficients of associated with the coefficients of $
f(\mathbf{x}).$
Therefore, in the following analysis, we only need to solve the following inhomogeneous polynomial optimization problem
\begin{equation}\label{e3}
\begin{aligned}
\min_{\x\in \mathbb{R}^n} &~f(\x) =\langle \mathcal{T}_f, M_d(\tilde{\mathbf{x}}\rangle=\mathcal{T}_f\tilde{\mathbf{x}}^4,\\
\textup{s.t.}   &~\|\textbf{\text{x}}\| = 1.
\end{aligned}
\end{equation}

\subsection{The corresponding multilinear optimization with spherical constraints}

In the following analysis, we denote the corresponding multilinear function with the same coefficient tensor of $f(\x)$ in (\ref{e3}) by
$$
F(\x,\y,\z,\w):=\langle\T_f, \tilde{\mathbf{x}}\circ\tilde{\mathbf{y}}\circ\tilde{\mathbf{z}}\circ\tilde{\mathbf{w}}\rangle= \mathcal{T}_f\tilde{\mathbf{x}}\tilde{\mathbf{y}}\tilde{\mathbf{z}}\tilde{\mathbf{w}},
$$
where $\tilde{\x}=(1, \x^\top)^\top,\tilde{\y}=(1, \y^\top)^\top,\tilde{\z}=(1, \z^\top)^\top,\tilde{\w}=(1, \w^\top)^\top\in\mathbb{R}^{n+1}$.
Therefore, the corresponding multilinear optimization problem can be denoted by
\begin{equation}\label{e4}
\begin{aligned}
\min &~
F(\x,\y,\z,\w):= \mathcal{T}_f\tilde{\mathbf{x}}\tilde{\mathbf{y}}\tilde{\mathbf{z}}\tilde{\mathbf{w}}, \\
\textup{s.t.}   &~ \|\mathbf{x}\|=\| \mathbf{y} \|=\| \mathbf{z} \|=\| \mathbf{w} \| = 1.
\end{aligned}
\end{equation}

The following theorem presents a sufficient condition to ensure that both ($\ref{e3}$) and ($\ref{e4}$) are equivalent with respect to the same optimal value. This is motivated by the Theorem 6 in \cite{Chen2021}, where it is proved that the homogeneous polynomial with degree four are equivalent with the corresponding multilinear polynomial optimization problem. To move on,
the following assumption is needed:
$$
{\bf(A)}: \mbox{The}~m\mbox{-th}~\mbox{order~polynomial}~f(\tilde{u})=\T_f\tilde{\u}^4~\mbox{is~concave~for~all}~\tilde{u}\in\mathbb{R}^{n+1}.
$$

\begin{theorem}\label{theorem1}
Suppose that the assumption {\bf(A)} holds. Let $\mathcal{T}_f \in\mathbf{S}^{4,n+1}$ be the tensor composed by the coefficients of $f(\x)$. Then ($\ref{e3}$) and ($\ref{e4}$) are equivalent such that
	\begin{equation}\label{e5}
	\min_{\|\mathbf{x}\|=1} f(\x)=\min_{\|\mathbf{x}\|,\|\mathbf{y}\|,\|\mathbf{z}\|,\|\mathbf{w}\|=1} F(\x, \y, \z, \w). \\
	\end{equation}
\end{theorem}
\proof To prove the conclusion, we first introduce the following bi-quadratic optimization problem:
\begin{equation*}
\begin{aligned}
\min &~ \bar{f}(\x,\y)=\langle\T_f,\tilde{\mathbf{x}}\circ\tilde{\mathbf{x}}\circ\tilde{\mathbf{y}}\circ\tilde{\mathbf{y}}\rangle=\mathcal{T}_f \tilde{\mathbf{x}} \tilde{\mathbf{x}} \tilde{\mathbf{y}} \tilde{\mathbf{y}}\\
\textup{s.t.}   &~  \|\mathbf{x}\| = \|\mathbf{y}\| =1,
\end{aligned}
\end{equation*}
where $\tilde{\mathbf{x}} = (1, \mathbf{x}^\top)^\top$, $\tilde{\mathbf{y}} = (1, \mathbf{y}^\top)^\top$.
First of all, it is obvious that
\begin{equation}\label{e6}
\min_{\|\mathbf{x}\|=1} f(\x) \geq
\min_{ \|\mathbf{x}\|,\|\mathbf{y}\| = 1 } \bar{f}(\x,\y)  \geq
 \min_{\|\mathbf{x}\|,\|\mathbf{y}\|,\|\mathbf{z}\|,\|\mathbf{w}\|=1} F(\x,\y,\z,\w).
\end{equation}

On the one hand, by assumption {\bf(A)}, it is equivalent with that for all $\tilde{\u}, \tilde{\v}\in\mathbb{R}^{n+1}$,
\begin{equation}\label{e7}
\tilde{\v}^\top\nabla^2f(\tilde{\u})\tilde{\v}=12\mathcal{T}_f\tilde{\mathbf{u}}\tilde{\mathbf{u}}\tilde{\mathbf{v}}\tilde{\mathbf{v}} \leq 0.
\end{equation}
For any $\mathbf{x}, \mathbf{y} \in \mathbb{R}^n$, denote $\tilde{\u}_1=(1,\x^\top)^\top, \tilde{\v}_2=(1,\y^\top)^\top$. Plugging $\tilde{\u}=\tilde{\u}_1+\tilde{\v}_2$, $\tilde{\v}=\tilde{\u}_1-\tilde{\v}_2$ into (\ref{e7}), it follows that
\begin{align*}
0 &\geq 	\langle\mathcal{T}_f, \tilde{\mathbf{u}}\circ\tilde{\mathbf{u}}\circ\tilde{\mathbf{v}}\circ\tilde{\mathbf{v}}\rangle \\
 &= \langle\mathcal{T}_f, (\tilde{\mathbf{u}}_1+\tilde{\mathbf{v}}_2) (\tilde{\mathbf{u}}_1+\tilde{\mathbf{v}}_2) (\tilde{\mathbf{u}}_1-\tilde{\mathbf{v}}_2) (\tilde{\mathbf{u}}_1-\tilde{\mathbf{v}}_2) \\
&= \langle\mathcal{T}_f, \tilde{\mathbf{u}}_1\circ\tilde{\mathbf{u}}_1\circ\tilde{\mathbf{u}}_1\circ\tilde{\mathbf{u}}_1\rangle + \langle\mathcal{T}_f, \tilde{\mathbf{v}}_2\circ \tilde{\mathbf{v}}_2\circ\tilde{\mathbf{v}}_2\circ\tilde{\mathbf{v}}_2\rangle-2\langle\mathcal{T}_f,\tilde{\mathbf{u}}_1\circ\tilde{\mathbf{u}}_1\circ\tilde{\mathbf{v}}_2\circ \tilde{\mathbf{v}}_2\rangle\\
&= \mathcal{T}_f \tilde{\mathbf{u}}_1^4 + \mathcal{T}_f \tilde{\mathbf{v}}_2^4 - 2\mathcal{T}_f \tilde{\mathbf{u}}_1^2\tilde{\mathbf{v}}_2^2,
\end{align*}
which implies that $2\mathcal{T}_f \tilde{\mathbf{u}}_1^2\tilde{\mathbf{v}}_2^2\geq \mathcal{T}_f \tilde{\mathbf{u}}_1^4 + \mathcal{T}_f \tilde{\mathbf{v}}_2^4\geq2\min\{\mathcal{T}_f \tilde{\mathbf{u}}_1^4,\mathcal{T}_f \tilde{\mathbf{v}}_2^4\}$ and
\begin{equation}\label{e8}
\min_{ \|\mathbf{x}\|,\|\mathbf{y}\| = 1 } \bar{f}(\x,\y)\geq\min_{\|\mathbf{x}\|=1} f(\x).
\end{equation}

On the other hand, for any $\tilde{\x}=(1,\x^\top)^\top, \tilde{\y}=(1,\y^\top)^\top,\tilde{\z}=(1,\z^\top)^\top, \tilde{\w}=(1,\w^\top)^\top\in \mathbb{R}^{n+1}$, by assumption {\bf(A)} again, we obtain that
\begin{align*}
0 &\geq  \mathcal{T}_f (\tilde{\mathbf{x}} +\tilde{\mathbf{y}})  (\tilde{\mathbf{x}} +\tilde{\mathbf{y}}) (\tilde{\mathbf{z}} -\tilde{\mathbf{w}}) (\tilde{\mathbf{z}} -\tilde{\mathbf{w}})  + \mathcal{T}_f (\tilde{\mathbf{x}} -\tilde{\mathbf{y}})  (\tilde{\mathbf{x}} -\tilde{\mathbf{y}}) (\tilde{\mathbf{z}} +\tilde{\mathbf{w}})  (\tilde{\mathbf{z}} +\tilde{\mathbf{w}}) \rangle \\
&= 2\mathcal{T}_f \tilde{\mathbf{x}} \tilde{\mathbf{x}}\tilde{\mathbf{z}} \tilde{\mathbf{z}}  +  2\mathcal{T}_f\tilde{\mathbf{x}} \tilde{\mathbf{x}} \tilde{\mathbf{w}} \tilde{\mathbf{w}} +  2\mathcal{T}_f\tilde{\mathbf{y}} \tilde{\mathbf{y}}\tilde{\mathbf{z}}\tilde{\mathbf{z}} +  2\mathcal{T}_f\tilde{\mathbf{y}}\tilde{\mathbf{y}}\tilde{\mathbf{w}}\tilde{\mathbf{w}}- 8\mathcal{T}_f\tilde{\mathbf{x}}\tilde{\mathbf{y}}\tilde{\mathbf{z}}\tilde{\mathbf{w}},
\end{align*}
which implies that
$$
8 \mathcal{T}_f\tilde{\mathbf{x}}\tilde{\mathbf{y}}\tilde{\mathbf{z}}\tilde{\mathbf{w}} \geq 8\min \left\{\mathcal{T}_f \tilde{\mathbf{x}} \tilde{\mathbf{x}}\tilde{\mathbf{z}} \tilde{\mathbf{z}} ,\mathcal{T}_f\tilde{\mathbf{x}} \tilde{\mathbf{x}} \tilde{\mathbf{w}} \tilde{\mathbf{w}},\mathcal{T}_f\tilde{\mathbf{y}} \tilde{\mathbf{y}}\tilde{\mathbf{z}}\tilde{\mathbf{z}},\mathcal{T}_f\tilde{\mathbf{y}}\tilde{\mathbf{y}}\tilde{\mathbf{w}}\tilde{\mathbf{w}} \right\}.
$$
Combining this with (\ref{e8}), and the arbitrariness of $\mathbf{x}$,$\mathbf{y}$,$\mathbf{z}$,$\mathbf{w}$, it follows that
$$
\min_{\|\mathbf{x}\|=1} f(\x) \leq
\min_{ \|\mathbf{x}\|,\|\mathbf{y}\| = 1 } \bar{f}(\x, \y)  \leq
\min_{\|\mathbf{x}\|,\|\mathbf{y}\|,\|\mathbf{z}\|,\|\mathbf{w}\|=1} F(\x, \y, \z, \w).
$$
Together with (\ref{e6}), it immediately leads to the assertion of this theorem. \qed

Generally speaking, the quartic function $f(\tilde{\u})=\T_f\tilde{\u}^4, \tilde{\u}\in\mathbb{R}^{n+1}$ is not necessarily concave for
any given symmetric tensor $\mathcal{T}_f$. To guarantee the concavity for given tenor $\mathcal{T}_f$, we introduce an augmented function, which includes an extra shift term :
$$
f_{\alpha}(\tilde{\mathbf{u}}) := \mathcal{T}_f\tilde{\mathbf{u}}^4 - \alpha \|\tilde{\mathbf{u}}\|^4.
$$
Then, the corresponding optimization problem of (\ref{e3}) takes the form:
\begin{equation}\label{e9}
\begin{aligned}
\min_{ \mathbf{x}  \in \mathbb{R}^n}~
f_{\alpha}(\x) := \mathcal{T}_f\tilde{\x}^4 - \alpha \|\tilde{\mathbf{x}}\|^4~~~
\textup{s.t.}~ \| \mathbf{x} \| = 1
\end{aligned}
\end{equation}
where $\alpha > 0$ can be regarded as a shift parameter to control the concavity of $f_{\alpha}(\tilde{\mathbf{x}})$.
Note that, the original (\ref{e3}) and the augmented problem (\ref{e9}) share the same optimal solution i.e.,
$$
\mathop{\arg\min}_{\|\mathbf{x}\|=1} f(\x) = \mathop{\arg\min}_{\|\mathbf{x}\|=1} f_\alpha(\x).
$$
To move on, the following theorem shows the concavity of the function $f_{\alpha}(\tilde{\x})$.
\begin{theorem}\label{theorem2}
For any given symmetric tensor $\mathcal{T}_f \in \mathbf{S}^{4,n+1}$. If $\alpha \geq \| \mathcal{T}_f \|$,
then the quartic function $f_{\alpha}(\tilde{\mathbf{u}})$ is concave with respect to $\tilde{\mathbf{u}} \in \mathbb{R}^{n+1}$.
\end{theorem}
\proof
For any $\tilde{\u}\in \mathbb{R}^{n+1}$, it is well-known that $f_{\alpha} (\tilde{\mathbf{u}} )$ is concave if and only if the matrix $ \nabla^2f_{\alpha}(\tilde{\mathbf{u}})$ is negative semi-definite i.e.,
$$
\tilde{\mathbf{v}}^\top \nabla^2f_{\alpha}(\tilde{\mathbf{u}})\tilde{\v}\leq 0, \quad \forall~ \tilde{\v}\in\mathbb{R}^{n+1}.
$$
By a direct computation, we have that
$$
\nabla^2f_{\alpha}(\tilde{\mathbf{u}})=12\mathcal{T}_f\tilde{\u}^2-8\alpha\tilde{\u}\tilde{\u}^\top-4\alpha\|\tilde{\u}\|^2I,
$$
where $I$ is the identity matrix with proper dimension. Then, for any $\tilde{\v}\in\mathbb{R}^{n+1}$, by Cauchy-Schwarz inequality, it follows that
$$
\begin{aligned}
\tilde{\mathbf{v}}^\top \nabla^2f_{\alpha}(\tilde{\mathbf{u}})\tilde{\v}&=12\mathcal{T}_f\tilde{\u}^2\tilde{\v}^2-8\alpha(\tilde{\u}^\top\tilde{\v})^2-4\alpha\|\tilde{\u}\|^2\|\tilde{\v}\|^2\\
&\leq (12\|\mathcal{T}_f\|-12\alpha)\|\tilde{\u}\|^2\|\tilde{\v}\|^2\leq 0,
\end{aligned}
$$
and the desired result follows.\qed

To propose the alternating algorithm, we need the corresponding multilinear optimization problem of (\ref{e9}).
To avoid introducing additionally redundant mathematical symbols, for any $\mathbf{x}$,$\mathbf{y}$,$\mathbf{z}$,$\mathbf{w}\in\mathbb{R}^n$, we denote the corresponding
multilinear function of $f_{\alpha}(\x)$ by $F_{\alpha}(\x,\y,\z,\w)$. Then, we have the following multilinear optimization problem:
\begin{equation}\label{e10}
\begin{aligned}
\min &\hspace{0.5cm}
F_\alpha(\x, \y, \z, \w) := \mathcal{T}_f\tilde{\mathbf{x}}\tilde{\mathbf{y}}\tilde{\mathbf{z}}\tilde{\mathbf{w}} - \alpha \langle \tilde{\mathbf{x}}, \tilde{\mathbf{y}}\rangle\langle \tilde{\mathbf{z}}, \tilde{\mathbf{w}}\rangle, \\
\textup{s.t.}   &\hspace{0.5cm} \| \mathbf{x} \| ,  \| \mathbf{y} \| ,  \| \mathbf{z} \| ,  \| \mathbf{w} \| = 1,
\end{aligned}
\end{equation}
where $\tilde{\x}=(1,\x^\top)^\top, \tilde{\y}=(1,\y^\top)^\top,\tilde{\z}=(1,\z^\top)^\top, \tilde{\w}=(1,\w^\top)^\top\in \mathbb{R}^{n+1}$.
As a direct result of Theorem \ref{theorem1} and Theorem \ref{theorem2}, we immediately have  the following theorem, which gives a sufficient condition to guarantee that both (\ref{e9}) and (\ref{e10}) have the same optimal value.
\begin{theorem}\label{theorem3}
	For $\mathcal{T}_f \in \mathbf{S}^{4,n+1}$, let $\alpha \geq \|\mathcal{T}_f\|$ be any positive constant, then we have the following results:
	\begin{equation}\label{e11}
	\min_{\|\mathbf{x}\|=1} f_{\alpha}(\x)=\min_{\|\mathbf{x}\|,\|\mathbf{y}\|,\|\mathbf{z}\|,\|\mathbf{w}\|=1} F_\alpha(\x, \y, \z, \w).
	\end{equation}
\end{theorem}

By Theorem \ref{theorem3} and the equivalence between (\ref{e3}) and (\ref{e9}), we immediately have that finding optimal solutions of (\ref{e3}) amounts to solving its augmented problem (\ref{e9}). Moreover, since the augmented multilinear (\ref{e10}) and augmented inhomogeneous form(\ref{e9}) have the same optimal value, we can accordingly
propose a proximal alternating minimization (PAM) algorithm inspired by exploiting the multi-block structure to get an approximate optimal value of the problem (\ref{e3}). For the sake of simplicity, define the following two sets:
$$
\mathbb{D} := \left\{\mathbf{x} \in \mathbb{R}^{n} ~|~ \|\mathbf{x}\| = 1\right\},~~
\mathbb{G} := \left\{\tilde{\mathbf{x}}={(1,{\mathbf{x}^\top)}}^\top \in \mathbb{R}^{n+1} ~|~ \| \mathbf{x} \|=1 \right\}.
$$
Next, we propose the following PAM algorithm.
\begin{algorithm}[!htbp]
	\caption{(Proximal Alternating Minimization Algorithm for \eqref{e3})}\label{alg1}
	\begin{algorithmic}[1]
		\STATE Let $\alpha \geq \|\mathcal{T}_f\|$, $\epsilon > 0$,
			$\tilde{\mathbf{x}}^{(0)}, \tilde{\mathbf{y}}^{(0)}, \tilde{\mathbf{z}}^{(0)}, \tilde{\mathbf{w}}^{(0)} \in \mathbb{G}$.
			 $\gamma_i \geq 0$ for $i = 1, 2, 3, 4$.
			\FOR{$k = 0, 1, 2, \ldots, n$}
		\STATE Update $(\tilde{\mathbf{x}}^{(k+1)}, \tilde{\mathbf{y}}^{(k+1)}, \tilde{\mathbf{z}}^{(k+1)}, \tilde{\mathbf{w}}^{(k+1)})$ sequentially via:
			$$
\left\{
			\begin{aligned}
			\tilde{\mathbf{x}}^{(k+1)} &= \arg\min_{{\tilde{\mathbf{x}}} \in \mathbb{G}}  F_\alpha(\mathbf{x}, \mathbf{y}^{(k)}, \mathbf{z}^{(k)}, \mathbf{w}^{(k)}) + \frac{\gamma_1}{2} \|\tilde{\mathbf{x}} - \tilde{\mathbf{x}}^{(k)}\|^2, \\
			\tilde{\mathbf{y}}^{(k+1)} &= \arg\min_{{\tilde{\mathbf{y}}} \in \mathbb{G}} F_\alpha(\mathbf{x}^{(k+1)}, \mathbf{y}, \mathbf{z}^{(k)}, \mathbf{w}^{(k)}) + \frac{\gamma_2}{2} \|\tilde{\mathbf{y}} - \tilde{\mathbf{y}}^{(k)}\|^2, \\
			\tilde{\mathbf{z}}^{(k+1)} &= \arg\min_{{\tilde{\mathbf{z}}} \in \mathbb{G}} F_\alpha(\mathbf{x}^{(k+1)},\mathbf{y}^{(k+1)}, \mathbf{z}, \mathbf{w}^{(k)}) + \frac{\gamma_3}{2} \|\tilde{\mathbf{z}} - \tilde{\mathbf{z}}^{(k)}\|^2, \\
			\tilde{\mathbf{w}}^{(k+1)} &= \arg\min_{{\tilde{\mathbf{w}}} \in \mathbb{G}} F_\alpha(\mathbf{x}^{(k+1)},\mathbf{y}^{(k+1)}, \mathbf{z}^{(k+1)}, \mathbf{w}) + \frac{\gamma_4}{2} \|\tilde{\mathbf{w}} - \tilde{\mathbf{w}}^{(k)}\|^2.
			\end{aligned}
\right.
			$$
		\STATE  Update $\tilde{\mathbf{u}}^{(k+1)}=\mathop{\arg\min} \left\{  f_\alpha(\tilde{\mathbf{x}}^{(k+1)}), f_\alpha (\tilde{\mathbf{y}}^{(k+1)}), f_\alpha (\tilde{\mathbf{z}}^{(k+1)}), f_\alpha(\tilde{\mathbf{w}}^{(k+1)}) \right\}  .$
		\STATE	If {$|f_\alpha(\mathbf{u}^{(k+1)})-f_\alpha(\mathbf{u}^{(k)})| \leq \epsilon$},
stop and return an approximate solutions $\u^{k+1}\in \mathbb{D}$.
			\ENDFOR
	\end{algorithmic}
\end{algorithm}

From Algorithm 1, two sequences $\{\tilde{\mathbf{t}}^{(k)} \} $ and $\{ \mathbf{t}^{(k)} \}$ will be generated, where
$\mathbf{t}^{(k)}=(\x^{(k)},\y^{(k)},\z^{(k)},\w^{(k)})$,
$$
\tilde{\mathbf{t}}^{(k)} =
(\tilde{\mathbf{x}}^{(k)}, \tilde{\mathbf{y}}^{(k)}, \tilde{\mathbf{z}}^{(k)}, \tilde{\mathbf{w}}^{(k)} ) =
 \left(\left( \begin{matrix} 1 \\ \mathbf{x}^{(k)} \end{matrix} \right) , \left( \begin{matrix} 1 \\ \mathbf{y}^{(k)} \end{matrix} \right) , \left( \begin{matrix} 1 \\ \mathbf{z}^{(k)} \end{matrix} \right) , \left( \begin{matrix} 1 \\ \mathbf{w}^{(k)} \end{matrix} \right) \right).
$$
Then $F_\alpha(\x^{(k)}, \y^{(k)}, \z^{(k)}, \w^{(k)})$ can be denoted by $\{ F_\alpha(\mathbf{t}^{(k)})\}$ for simple.
We now present some convergence properties for Algorithm \ref{alg1}.
\begin{theorem}\label{thm4}
	Let $\{ \tilde{\mathbf{t}}^{(k)}\}$,$\{\mathbf{t}^{(k)}\}$ and $\{ \tilde{\mathbf{u}}^{(k)}\}$be sequences generated by Algorithm \ref{alg1} with any initial point 	$\tilde{\mathbf{x}}^{(0)}, \tilde{\mathbf{y}}^{(0)}, \tilde{\mathbf{z}}^{(0)}, \tilde{\mathbf{w}}^{(0)} \in \mathbb{G}$. Then, we have the following results.
	\begin{itemize}
	\itemindent 1em
	\item[{\rm (i)}]
		The sequence $\{ F_\alpha(\mathbf{t}^{(k)})\}$ is monotonically decreasing, and there is a real number $f^*_\alpha$ such that
		$$
		\lim_{k \to \infty} f_\alpha(\mathbf{t}^{(k)}) = f^*_\alpha.
		$$
    \item[{\rm (ii)}]
		It holds that $\|  \mathbf{t}^{(k+1)} - \mathbf{t}^{(k)}  \|  \to 0$ as $k \to \infty$.
	\item[{\rm\ (iii)}]
		There is a subsequence $\{ \tilde{\mathbf{u}}^{(p_k)} \} \subseteq \{ \tilde{\mathbf{u}}^{(k)} \}$ converges to a real number $f^*$:
		\[
		\lim_{k \to \infty} f_\alpha(\mathbf{u}^{(p_k)}) = f^*.
		\]
	\item[{\rm (iv)}]	
		The sequence $\{\tilde{\mathbf{u}}^{(k)} \}$ has limit points, and each of its convergent subsequence converges to $\tilde{\mathbf{u}}^* =\left( \begin{matrix} 1 \\ \mathbf{u}^* \end{matrix} \right) $ such that \[
		f_\alpha(\mathbf{u}^*) \leq f^*_\alpha.
		\]
	\end{itemize}
\end{theorem}
\proof
	Let $\gamma_{\text{min}} := \min\{\gamma_1,\gamma_2,\gamma_3,\gamma_4\} > 0$, it follows from the iterative schemes (i.e., step 3) of Algorithm \ref{alg1} that
\begin{equation}\label{eq:4}
\begin{aligned}
&F_\alpha(\mathbf{t}^{(k)}) = F_\alpha(\mathbf{x}^{(k)}, \mathbf{y}^{(k)}, \mathbf{z}^{(k)}, \mathbf{w}^{(k)}) \\
	&\geq F_\alpha(\mathbf{x}^{(k+1)},\mathbf{y}^{(k)}, \mathbf{z}^{(k)},\mathbf{w}^{(k)}) + \frac{\gamma_1}{2} \|  \tilde{\mathbf{x}}^{(k+1)} - \tilde{\mathbf{x}}^{(k)} \| ^2 \\
	&\geq F_\alpha(\mathbf{x}^{(k+1)}, \mathbf{y}^{(k+1)}, \mathbf{z}^{(k)}, \mathbf{w}^{(k)}) + \frac{\gamma_1}{2} \|  \tilde{\mathbf{x}}^{(k+1)} - \tilde{\mathbf{x}}^{(k)}  \| ^2 \\
&~~~+\frac{\gamma_2}{2} \|  \tilde{\mathbf{y}}^{(k+1)} - \tilde{\mathbf{y}}^{(k)}  \| ^2 \geq F_\alpha(\mathbf{x}^{(k+1)}, \mathbf{y}^{(k+1)}, \mathbf{z}^{(k+1)}, \mathbf{w}^{(k)})\\
&~~~ + \frac{\gamma_1}{2} \|  \tilde{\mathbf{x}}^{(k+1)} - \tilde{\mathbf{x}}^{(k)}  \| ^2+ \frac{\gamma_2}{2} \|  \tilde{\mathbf{y}}^{(k+1)} - \tilde{\mathbf{y}}^{(k)}  \| ^2 + \frac{\gamma_3}{2} \|  \tilde{\mathbf{z}}^{(k+1)} - \tilde{\mathbf{z}}^{(k)}  \| ^2 \\
	&\geq F_\alpha(\mathbf{x}^{(k+1)}, \mathbf{y}^{(k+1)}, \mathbf{z}^{(k+1)}, \mathbf{w}^{(k+1)}) + \frac{\gamma_{\text{min}}}{2} \|  \tilde{\mathbf{t}}^{(k+1)} - \tilde{\mathbf{t}}^{(k)}  \| ^2 \\
	&= F_\alpha(\mathbf{t}^{(k+1)})+\frac{\gamma_{\text{min}}}{2} \|  \tilde{\mathbf{t}}^{(k+1)} -\tilde{\mathbf{t}}^{(k)}  \| ^2\geq F_\alpha(\mathbf{t}^{(k+1)}) ,
\end{aligned}
\end{equation}
which shows that the sequence $\{F_\alpha(\mathbf{t}^{(k)})\}$ is monotonically decreasing. On the other hand, by the fact that the function $F_\alpha(\mathbf{x}, \mathbf{y}, \mathbf{z}, \mathbf{w})$ is continuous and $\mathbb{G}$ is compact, we know that the sequence $\{ F_\alpha(\mathbf{t}^{(k)})\}$ is bounded, which implies that there is a unique limit $f^*_\alpha$ such that
	$$
	\lim_{k \to \infty} F_\alpha(\mathbf{t}^{(k)}) = f^*_\alpha.
	$$
	
To prove (ii), by (\ref{eq:4}), we obtain that
$$
\frac{\gamma_{\text{min}}}{2} \|  \mathbf{t}^{(k+1)} - \mathbf{t}^{(k)}  \| ^2 = \frac{\gamma_{\text{min}}}{2} \|  \tilde{\mathbf{t}}^{(k+1)} - \tilde{\mathbf{t}}^{(k)}  \| ^2 \leq F_\alpha(\mathbf{t}^{(k)}) - F_\alpha(\mathbf{t}^{(k+1)}).
$$
Combining this with the fact that $\{ F_\alpha(\mathbf{t}^{(k)})\}$ is bounded and monotonically decreasing, we immediately arrive at the conclusion that
	$$
	\|  \mathbf{t}^{(k+1)} - \mathbf{t}^{(k)}  \|  \to 0 \text{ as } k \to \infty.
	$$

By the results of (i), (ii) and Algorithm \ref{alg1}, the results of (iii), (iv) can be naturally proved and the desired results hold.\qed

If an infinite sequence is generated in Algorithm \ref{alg1}, we have the following result.
   \begin{theorem} \label{thm5}
	Let $\{\mathbf{t}^{(k)} \mid \mathbf{t}^{(k)} = (\mathbf{x}^{(k)}, \mathbf{y}^{(k)}, \mathbf{z}^{(k)}, \mathbf{w}^{(k)}) \}$  be the infinite sequence generated by Algorithm \ref{alg1}. Then, any limit point $\mathbf{t}^* = (\mathbf{x}^*, \mathbf{y}^*, \mathbf{z}^*, \mathbf{w}^*) $ of $\{\mathbf{t}^{(k)}\}$ is a stationary point of the augmented multilinear problem $(\ref{e10})$.
	\end{theorem}
\proof
		Due to the compactness of $\mathbb{G}$, it is not difficult to see that the sequence $\{\mathbf{t}^{(k)}\}$ generated by Algorithm \ref{alg1} is bounded, which further implies that such a sequence has at least one limit point. Suppose that $\{\mathbf{t}^{(p_k)}\}$ is a subsequence of $\{\mathbf{t}^{(k)}\}$  with limit point $\mathbf{t}^* = (\mathbf{x}^*, \mathbf{y}^*,\mathbf{z}^*, \mathbf{w}^*) $, i.e.,
		\begin{align*}
    	\tilde{\mathbf{t}}^{(p_k)} &= \left( \left( \begin{matrix} 1 \\ \mathbf{x}^{(p_k)} \end{matrix} \right) , \left( \begin{matrix} 1 \\ \mathbf{y}^{(p_k)} \end{matrix} \right) , \left( \begin{matrix} 1 \\ \mathbf{z}^{(p_k)} \end{matrix} \right) , \left( \begin{matrix} 1 \\ \mathbf{w}^{(p_k)} \end{matrix} \right) \right) \\
	&\to
	  \left(\left( \begin{matrix} 1 \\ \mathbf{x}^* \end{matrix} \right) , \left( \begin{matrix} 1 \\ \mathbf{y}^* \end{matrix} \right) , \left( \begin{matrix} 1 \\ \mathbf{z}^* \end{matrix} \right) , \left( \begin{matrix} 1 \\ \mathbf{w}^* \end{matrix} \right) \right) \text{ as } k \to \infty.
	  \end{align*}
		It follows from (ii) of Theorem \ref{thm4} that
		$$
		\| \mathbf{t}^{(p_k+1)} -\mathbf{t}^*  \| \leq	\| \mathbf{t}^{(p_k+1)} - \mathbf{t}^{(p_k)} \| + \| \mathbf{t}^{(p_k)} - \mathbf{t}^*  \|
		\to 0 \text{ as } k \to \infty.
		$$		
		
		To verify that $\mathbf{t}^*$ is a stationary point of (\ref{e10}), we just need to prove the following variational inequality holds:
		\begin{equation}\label{eq:5}
		\left\langle  \nabla F_\alpha(\mathbf{t}^*), \mathbf{t}- \mathbf{t}^* \right\rangle  \geq 0, \quad \forall~ \mathbf{t} \in \mathbb{D} \times \mathbb{D}  \times \mathbb{D}  \times \mathbb{D} .
		\end{equation}
First of all, for the $\mathbf{x}$-subproblem of Algorithm \ref{alg1}, by (\ref{eq:4}), we have that
		$$
		F_\alpha( \mathbf{x}, \mathbf{y}^{(p_k)}, \mathbf{z}^{(p_k)}, \mathbf{w}^{(p_k)}  ) \geq F_\alpha(\mathbf{x}^{(p_{k+1})},\mathbf{y}^{(p_k)}, \mathbf{z}^{(p_k)}, \tilde{\mathbf{w}}^{(p_k)})\geq F_{\alpha}(\t^*) , \quad \forall~ \tilde{\mathbf{x}} \in \mathbb{G}.
		$$
		Note that
		$$
		\nabla_{ \mathbf{x} } F_\alpha(\mathbf{t}^{(k)})=\left( \begin{matrix} \mathbf{0} ,& \mathbf{I_n} \end{matrix} \right) \mathcal{T}_f \tilde{\mathbf{y}}^{(k)} \tilde{\mathbf{z}}^{(k)} \tilde{\mathbf{w}}^{(k)} -\alpha \langle  \tilde{\mathbf{z}}^{(k)}, \tilde{\mathbf{w}}^{(k)} \rangle  	\left( \begin{matrix} \mathbf{0} ,& \mathbf{I_n} \end{matrix} \right)  \tilde{\mathbf{y}}^{(k)},
		$$
where $\0\in\mathbb{R}^n$ and $I_n$ is the identity matrix with dimension $n\times n$.
		Consequently, let $k \rightarrow \infty$ and it follows that
		\begin{align*}
		\left\langle  \nabla_{\mathbf{x} } F_\alpha(\mathbf{t}^*), \mathbf{x} - \mathbf{x}^* \right\rangle
		&=\left\langle  ( \begin{matrix} \mathbf{0} ,& \mathbf{I_n} \end{matrix} ) \mathcal{T}_f \tilde{\mathbf{y}}^* \tilde{\mathbf{z}}^* \tilde{\mathbf{w}}^* ,\mathbf{x}- \mathbf{x}^* \right\rangle  -\alpha \langle \tilde{\mathbf{z}}^*, \tilde{\mathbf{w}}^* \rangle 	\left\langle  ( \begin{matrix} \mathbf{0} ,& \mathbf{I_n} \end{matrix} )  \tilde{\mathbf{y}}^*, \mathbf{x}-\mathbf{x}^* \right\rangle \\
		&=\mathcal{T}_f   \left( \begin{matrix} 0 \\ \mathbf{x-x^*} \end{matrix} \right)   \left( \begin{matrix} 1 \\ \mathbf{y}^* \end{matrix} \right)  \left( \begin{matrix} 1 \\ \mathbf{z}^* \end{matrix} \right)  \left( \begin{matrix} 1 \\ \mathbf{w}^* \end{matrix} \right)   -\alpha \langle \tilde{\mathbf{z}}^*, \tilde{\mathbf{w}}^* \rangle \langle \mathbf{y}^*, \mathbf{x}-\mathbf{x}^* \rangle\\
		&=\mathcal{T}_f (\tilde{\mathbf{x}}-\tilde{\mathbf{x}}^*)  \tilde{\mathbf{y}}^* \tilde{\mathbf{z}}^* \tilde{\mathbf{w}}^* -\alpha \langle \tilde{\mathbf{z}}^*, \tilde{\mathbf{w}}^* \rangle \langle \tilde{\mathbf{x}}-\tilde{\mathbf{x}}^*,\tilde{\mathbf{y}}^* \rangle\\
		&= F_{\alpha}(\mathbf{x}-\mathbf{x}^*, \mathbf{y}^*, \mathbf{z}^*, \mathbf{w}^*)\\
		&= F_{\alpha}(\mathbf{x}, \mathbf{y}^*, \mathbf{z}^*, \mathbf{w}^*) - F_{\alpha}(\t^*)\geq 0.
		\end{align*}
From a similar proof with above, it follows that, for all $\y, \z, \w\in\mathbb{D}$,
$$
\left\langle  \nabla_{\mathbf{y} } F_\alpha(\mathbf{t}^*), \mathbf{y} - \mathbf{y}^* \right\rangle\geq 0,~
	    	\left\langle  \nabla_{\mathbf{z} } F_\alpha(\mathbf{t}^*), \mathbf{z} - \mathbf{z}^* \right\rangle    \geq 0, ~
			\left\langle  \nabla_{\mathbf{w} } f_\alpha(\mathbf{t}^*), \mathbf{w} - \mathbf{w}^* \right\rangle    \geq 0.
$$
Invoking the fact that $
		 \nabla F_\alpha(\mathbf{t}^*) = ( \nabla_{\mathbf{x}} F_\alpha(\mathbf{t}^*), \nabla_{\mathbf{y}} F_\alpha(\mathbf{t}^*), \nabla_{\mathbf{z}} F_\alpha(\mathbf{t}^*), \nabla_{\mathbf{w}} F_\alpha(\mathbf{t}^*) ),
		$
then we know that (\ref{eq:5}) holds and $\mathbf{t}^*$ is a stationary point of (\ref{e10}).
	\qed	

Since the linear independent constraint qualification holds automatically under spherical constraints, we immediately know that the optimal solution is a KKT point for each subproblem of Algorithm 1. Recalling the $\mathbf{x}$-subproblem of Algorithm 1, i.e.,
	$$
	\min_{{\mathbf{x}} \in \mathbb{D}}  F_\alpha(\mathbf{x}, \mathbf{y}^{(k)}, \mathbf{z}^{(k)}, \mathbf{w}^{(k)}) + \frac{\gamma_1}{2} \|\tilde{\mathbf{x}} - \tilde{\mathbf{x}}^{(k)}\|^2,
	$$	
we immediately have a pair of KKT points $(\mathbf{x}^*, \lambda^* )$ satisfying $\|\x^*\| = 1$ and
\begin{equation*}
\begin{aligned}
&\left( \begin{matrix} \mathbf{0} ,& \mathbf{I_n} \end{matrix} \right) \mathcal{T} _f \tilde{\mathbf{y}}^{(k)} \tilde{\mathbf{z}}^{(k)} \tilde{\mathbf{w}}^{(k)} -\alpha \langle  \tilde{\mathbf{z}}^{(k)}, \tilde{\mathbf{w}}^{(k)} \rangle  	\left( \begin{matrix} \mathbf{0} ,& \mathbf{I_n} \end{matrix} \right)  \tilde{\mathbf{y}}^{(k)}\\
&+\gamma_1\left( \begin{matrix} \mathbf{0} ,& \mathbf{I_n} \end{matrix} \right)(\tilde{\mathbf{x}}^*-\tilde{\mathbf{x}}^{(k)})-2\lambda^*\left( \begin{matrix} \mathbf{0} ,& \mathbf{I_n} \end{matrix} \right)\tilde{\mathbf{x}}^*=\0
\end{aligned}
\end{equation*}
By a direct computation, $\mathbf{x}^*$ can be expressed explicitly by
\begin{equation}\label{e14}
\begin{aligned}
\mathbf{x}^* &= \pm \frac{\left( \begin{matrix} \mathbf{0} ,& \mathbf{I_n} \end{matrix} \right) \mathcal{T}_f\tilde{\mathbf{y}}^{(k)} \tilde{\mathbf{z}}^{(k)} \tilde{\mathbf{w}}^{(k)} -\alpha \langle  \tilde{\mathbf{z}}^{(k)}, \tilde{\mathbf{w}}^{(k)} \rangle  \mathbf{y}^{(k)} - \gamma_1 \mathbf{x}^{(k)} }
{\| \left( \begin{matrix} \mathbf{0} ,& \mathbf{I_n} \end{matrix} \right) \mathcal{T}_f\tilde{\mathbf{y}}^{(k)} \tilde{\mathbf{z}}^{(k)} \tilde{\mathbf{w}}^{(k)} -\alpha \langle  \tilde{\mathbf{z}}^{(k)}, \tilde{\mathbf{w}}^{(k)} \rangle  \mathbf{y}^{(k)} - \gamma_1 \mathbf{x}^{(k)} \|}.
\end{aligned}
\end{equation}
Note that if the denominators of (\ref{e14}) equal zero, we may set \(\x^{(k+1)} = \x^{(k)}\) and adjust the parameters \(\gamma_1\) and \(\alpha\). Similarly, the KKT points for the subproblems with blocks $\y, \z, \w$ can be derived analogously. Since all subproblems in Algorithm \ref{alg1} admit finite analytic solutions, the algorithm is computationally efficient for spherically constrained problems.

Although an approximate optimal value for \eqref{e9} can be achieved by Algorithm \ref{alg1}, a solution produced by Algorithm \ref{alg1} may not be a KKT point of (\ref{e9}).
For this purpose, we consider a combination of BIM \cite{Chen2021-29} and Algorithm \ref{alg1}. Specifically, we first use Algorithm \ref{alg1} to produce an approximate optimal solution $\mathbf{u}^*$, which is the starting point of Algorithm \ref{alg2}.
\begin{algorithm}[!htbp]
	\caption{(Enhanced BIM for \eqref{e9})}\label{alg2}
	\begin{algorithmic}[1]
		\STATE Let $\alpha \geq \|\mathcal{T}_f \|$, $\epsilon > 0$, Take $\mathbf{x}^{(0)}=\mathbf{u}^*$, where $\mathbf{u}^*$ is obtained from Algorithm \ref{alg1}.
			\FOR{$k = 0, 1, 2, \ldots, N$}
		\STATE Update $\mathbf{x}^{(k+1)}$ via:
			$$
			\mathbf{x}^{(k+1)} =
			\begin{cases}
				\hspace{0.7cm} \mathbf{x}^{(k)}\hspace{0.6cm}, & \text{if  } \nabla_{\mathbf{x}} f_\alpha(\mathbf{x}^{(k)}) = \0, \\
				-\frac{\nabla_{\mathbf{x}} f_\alpha(\mathbf{x}^{(k)})}{\| \nabla_{\mathbf{x}} f_\alpha(\mathbf{x}^{(k)}) \|}\hspace{0.1cm}, & \text{otherwise}.
			\end{cases}
		$$
		\STATE Update $\tilde{\mathbf{x}}^{(k+1)}=\left( \begin{matrix} 1 \\ \mathbf{x}^{(k+1)} \end{matrix} \right). $
		\STATE If {$|f_\alpha(\mathbf{x}^{(k+1)}) - f_\alpha(\mathbf{x}^{(k)})| \leq \epsilon$}, stop and return a stationary point $\mathbf{x}^{(k+1)}$.
			\ENDFOR
	\end{algorithmic}
\end{algorithm}

\begin{theorem}\label{thm6}
 	Let $\{\tilde{\mathbf{x}}^{(k)}\}=\left\{ \left( \begin{matrix} 1 \\ \mathbf{x}^{(k)} \end{matrix} \right) \right\}$ be the sequence generated in Algorithm \ref{alg2}.  If $\alpha \geq \|\mathcal{T}_f \|$, then the following results hold.
	\begin{itemize} 	
		\itemindent 1em
 		\item[{\rm (i)}]
 		The sequence $ \{f_{\alpha}(\mathbf{x}^{(k)})\}$ is strictly decreasing.
 		\item[{\rm (ii)}]
 		Any cluster point of $\{\mathbf{x}^{(k)}\}$ is a KKT point of \eqref{e9}.
 		\item[{\rm (iii)}]
 		 If $f_{\alpha}(\mathbf{x})$ is strictly concave, then $\|\mathbf{x}^{(k+1)}-\mathbf{x}^{(k)}\| \to 0$ as $k \to \infty$.
	\end{itemize}
\end{theorem}
\proof Since $\alpha \geq \|\mathcal{T}_f \|$, then $f_ \alpha(\mathbf{x})$ is concave for all $\x\in\mathbb{R}^n$ and it holds that, for $k\in\mathbb{N}$,
 	\begin{equation}\label{eq:6}
 	f_\alpha(\mathbf{x}^{(k+1)}) - f_\alpha(\mathbf{x}^{(k)}) \leq \langle  \nabla_{\mathbf{x}} f_\alpha(\mathbf{x}^{(k)}), \mathbf{x}^{(k+1)} - \mathbf{x}^{(k)} \rangle  .
 	\end{equation}
If $\nabla_{\mathbf{x}} f_\alpha(\mathbf{x}^{(k)})=0$, then, it follows from Algorithm \ref{alg2} that $\mathbf{x}^{(k+1)} =\mathbf{x}^{(k)} $ and (\ref{eq:6}) that
$f_\alpha(\mathbf{x}^{(k+1)})\leq f_\alpha(\mathbf{x}^{(k)})$. Otherwise, i.e., $ \nabla_{\mathbf{x}} f_\alpha(\mathbf{x}^{(k)}) \neq 0 $, by (\ref{eq:6}) again, it follows that
 	\begin{equation*}
 	\begin{aligned}
     &f_\alpha(\mathbf{x}^{(k+1)}) - f_\alpha(\mathbf{x}^{(k)}) \leq \langle  \nabla_{\mathbf{x}} f_\alpha(\mathbf{x}^{(k)}), \mathbf{x}^{(k+1)} - \mathbf{x}^{(k)}  \rangle\\
      &=\left\langle  \nabla_{\mathbf{x}} f_\alpha(\mathbf{x}^{(k)}), 	-\frac{\nabla_{\mathbf{x}} f_\alpha(\mathbf{x}^{(k)})}{\| \nabla_{\mathbf{x}} f_\alpha(\mathbf{x}^{(k)}) \|} -\mathbf{x}^{(k)} \right\rangle \\&=
     -\|\nabla_{\mathbf{x}} f_\alpha (\mathbf{x}^{(k)}) \| - \langle  \nabla_{\mathbf{x}} f_\alpha(\mathbf{x}^{(k)}), \mathbf{x}^{(k)} \rangle\leq
     \|\nabla_{\mathbf{x}} f_\alpha (\mathbf{x}^{(k)}) \| ( \|\mathbf{x}^{(k)} \|-1 )=0,
 		\end{aligned}
 	\end{equation*}
 where the last inequality follows from the Cauchy-schwarz inequality. Combining this with the fact that $\|\x^{(k+1)}\|=\|\x^{(k)}\|=1$, we have that
\begin{equation}\label{eq:7}
f_\alpha(\mathbf{x}^{(k+1)}) - f_\alpha(\mathbf{x}^{(k)}) < 0 \Leftrightarrow \mathbf{x}^{(k+1)} \neq \mathbf{x}^{(k)},
\end{equation}
which means that $\{f_{\alpha}(\mathbf{x}^{(k)})\}$ is strictly decreasing.

To prove (ii), by the monotonicity of $\{ f_\alpha(\mathbf{x}^{(k)})\}$ and the compactness of $\mathbb{G}$, we know that $\{ f_\alpha(\mathbf{x}^{(k)})\}$ is convergent. By (\ref{eq:6}), we obtain that
 	\begin{equation}\label{eq:8}
 	\lim_{k \to \infty} \langle  \nabla_{ \mathbf{x}}f_\alpha(\mathbf{x}^{(k)}) , \mathbf{x}^{(k+1)} - \mathbf{x}^{(k)} \rangle  =0.
 	\end{equation}
 	Without loss of generality, assume $\tilde{{\mathbf{x}}}^*= \left( \begin{matrix} 1 \\ \mathbf{x}^* \end{matrix} \right) $ is an accumulation point of the subsequence $\{\tilde{\mathbf{x}}^{(k_j)}\}= \left \{ \left( \begin{matrix} 1 \\ \mathbf{x}^{(k_j)} \end{matrix} \right) \right\}$. Denote
$ L(\mathbf{x},\lambda) = \mathcal{T}_f  \tilde{\mathbf{x}}^4 - \lambda (\|\mathbf{x}\|^2-1).$ Then, it is enough to verify that there exist $\lambda^*\in\mathbb{R}$ such that $\nabla_{\mathbf{x}} L({\mathbf{x}}^*,\lambda^*)=\0 $, or equivalently,
 	$$
 		\left( \begin{matrix}\0, & \mathbf{I_n} \end{matrix} \right) \mathcal{T}_f\tilde{\x}^*\tilde{\x}^*\tilde{\x}^*= 2\lambda^* \mathbf{x}^*.
 		$$
We prove the conclusion from the following two cases.
 	If $\lim_{k \to \infty}  \nabla_{ \mathbf{x}}f_\alpha(\mathbf{x}^{(k_j)}) =\0$, it then from the continuity of $f_\alpha(\mathbf{x})$
 	that $\nabla_{\mathbf{x}}f_\alpha(\mathbf{x}^*)=\0,$ which means that
 	 	 \begin{equation}\label{eq:9}
 	 	\left( \begin{matrix} \0 ,& \mathbf{I_n} \end{matrix} \right)\mathcal{T}_f\tilde{\x}^*\tilde{\x}^*\tilde{\x}^*= \alpha \left( \begin{matrix} \0 ,& \mathbf{I_n} \end{matrix} \right) \tilde{\mathbf{x}}^* = \alpha \mathbf{x}^*.
 	 	 \end{equation}
 On the other hand, if $\lim_{k \to \infty}  \nabla_{\mathbf{x}}f_\alpha(\mathbf{x}^{(k_j)}) \neq \0$, by (\ref{eq:8}) again, we have that
 	 		\begin{equation*}
 	 	\begin{aligned}
 	 	\lim_{k \to \infty} \langle  \nabla_{\mathbf{x}} f_\alpha(\mathbf{x}^{(k_j)}), \mathbf{x}^{(k_j+1)} - \mathbf{x}^{(k_j)} \rangle &=\lim_{k \to \infty} \left\langle  \nabla_{\mathbf{x}} f_\alpha(\mathbf{x}^{(k_j)}), 	-\frac{\nabla_{\mathbf{x}} f_\alpha(\mathbf{x}^{(k_j)})}{\| \nabla_{\mathbf{x}} f_\alpha(\mathbf{x}^{(k_j)}) \|} - \mathbf{x}^{(k_j)} \right\rangle \\&=
 	 	-\|\nabla_{\mathbf{x}} f_\alpha(\mathbf{x}^*) \| - \langle  \nabla_{\mathbf{x}} f_\alpha(\mathbf{x}^*), \mathbf{x}^* \rangle  \\&=0,
 	 	\end{aligned}
 	 	\end{equation*}
which implies that
 	 	 \begin{equation}\label{eq:10}
 	 	 \nabla_{\mathbf{x}} f_\alpha(\mathbf{x}^*)=-\|\nabla_{\mathbf{x}} f_\alpha (\mathbf{x}^*) \| \mathbf{x}^*.
 	 	 \end{equation}
Plugging  $\nabla_{\mathbf{x}} f_\alpha(\mathbf{x}^*)=4\left( \begin{matrix} \0 ,& \mathbf{I_n} \end{matrix} \right)\mathcal{T}_f\tilde{\x}^*\tilde{\x}^*\tilde{\x}^*-8\alpha \mathbf{x}^*$ into (\ref{eq:10}) gives
 	 	\begin{equation*} 	 	
 	 	\left( \begin{matrix} \0 ,& \mathbf{I_n} \end{matrix} \right) \mathcal{T}_f\tilde{\x}^*\tilde{\x}^*\tilde{\x}^*= \frac{1}{4} (8\alpha-\|\nabla_{\mathbf{x}} f_\alpha(\mathbf{x}^*)\|)\mathbf{x}^*.
 	 	\end{equation*}
Combining (\ref{eq:9}) with (\ref{eq:10}), we know that $\mathbf{x}^*$ is a KKT point of \eqref{e9}.
 	 	
To prove (iii), if $f_\alpha(\mathbf{x})$ is strictly concave, we have that $f_\alpha(\mathbf{x}^{(k+1)})< f_\alpha(\mathbf{x}^{(k)})$ for any $\tilde{\mathbf{x}}^{(k)} \in \mathbb{G}$. 	
 	 	By (\ref{eq:7}) and the definition of $\mathbf{x}^{(k+1)}$(see step 3 of the Algorithm 2), we have
 	 	\begin{equation*}
 	 	\begin{aligned}
 	 	\| \mathbf{x}^{(k+1)} \|^2- \langle  \mathbf{x}^{(k+1)}, \mathbf{x}^{(k)} \rangle   &=
 	 	\langle  \mathbf{x}^{(k+1)},  \mathbf{x}^{(k+1)}-\mathbf{x}^{(k)} \rangle  \\&=-\frac{1}{\| \nabla_{\mathbf{x}} f_\alpha(\mathbf{x}^{(k)}) \|}
 	 	\langle  \nabla_{\mathbf{x}} f_\alpha(\mathbf{x}^{(k)}), 	\mathbf{x}^{(k+1)}- \mathbf{x}^{(k)} \rangle .
 	 	\end{aligned}
 	 	\end{equation*}
 	 	It then follows from (\ref{eq:8}) that
 	 	 	$$
 	 	 	\lim_{k \to \infty} \langle  \mathbf{x}^{(k+1)}, \mathbf{x}^{(k)} \rangle   = \lim_{k \to \infty} \| \mathbf{x}^{(k+1)} \|^2 = 1.
 	 	 	$$
 	 	 	Then, we conclude that
 	 	 	$$
 	 	 	\lim_{k \to \infty}
 	 	 	\| \mathbf{x}^{(k+1)}-\mathbf{x}^{(k)} \| = 0,
 	 	 	$$
and the desired results hold. \qed

\section{Numerical examples}\label{sec03}
In this section, we provide many numerical results to verify the efficacy of the proposed method on non-rotating BEC problems.
For this purpose, recall the references \cite{2022-1,2022-3,2022-28}, where the corresponding energy function is defined as below:
 		\begin{equation}\label{eq:Ex}
 	E(\phi(\mathbf{x})) := \int_{\mathbb{R}^d} \left[ \frac{1}{2} |\nabla \phi(\mathbf{x})|^2 + V(\mathbf{x})|\phi(x)|^2 + \frac{\beta}{2} |\phi(\mathbf{x})|^4 \right] \, {\mbox d}\mathbf{x}.
 	\end{equation}
Then the ground state of non-rotating BEC is usually defined as the minimizer of the following nonconvex minimization problem:
	\begin{equation}\label{eq:BEC}
 	\begin{aligned}
 	\min  &\hspace{0.5cm}   \int_{\mathbb{R}^d} \left[ \frac{1}{2} |\nabla \phi(\mathbf{x})|^2 + V(\mathbf{x})|\phi(x)|^2 + \frac{\beta}{2} |\phi(\mathbf{x})|^4 \right] \, {\mbox d}\mathbf{x}\\
 	\textup{s.t.}   &\hspace{0.5cm} \int_{\mathbb{R}^d} |\phi(\mathbf{x})|^2 \, d\mathbf{x} = 1, \quad E(\phi(\x)) < \infty
 	\end{aligned}
\end{equation}
 	where $\mathbf{x} \in \mathbb{R}^n$ is the spatial coordinate vector, $V(\mathbf{x})$ is an external trapping potential, and the given constant $\beta \in \mathbb{R}$ is the dimensionless interaction coefficient, see \cite{2022-3}. The optimal solution $\phi^*(\mathbf{x})$ is defined as the ground state. We only consider $\beta > 0$ in the following examples. For practical applications of BEC, the harmonic potential is often used as below \cite{2022-6,2022-7}:
 	$$
 	V(\mathbf{x}) =  \frac{1}{2}
 	\begin{cases}
 	\gamma^2 x^2, & d=1, \\
    \gamma_1^2 x_1^2 + \gamma_2^2 x_2^2, & d=2,
 	\end{cases}
 	$$
    where $\gamma, \gamma_1, \gamma_2$ are three given positive constants. Unless otherwise specified, let $\gamma = \gamma_1  = \gamma_2=1$ in the following analysis.

Due to the external trapping potential, the ground state of $(\ref{eq:BEC})$ exhibits exponential decay as $|\mathbf{x}| \to \infty$.
This behavior allows the energy function to be truncated from the entire space $\mathbb{R}^d$ to a sufficiently large bounded computational domain $\mathbb{D}$, such that the truncation error becomes negligible under Dirichlet boundary conditions \cite{2022-1,2022-3,2022-4}.
The finite difference method is then applied to discretize the energy function (\ref{eq:Ex}) and the constraints in (\ref{eq:BEC}) \cite{2022-4}. As a result, the Bose–Einstein condensate problem can be reduced to a finite-dimensional minimization problem with a spherical constraint.
Appropriate discretization schemes enable the numerical solution of (\ref{eq:BEC}). Several commonly used methods for computing the ground state include finite difference discretization \cite{BEC-21}, \cite{BEC-22} and pseudospectral approximation based on fast Fourier transform or discrete sinusoidal transform \cite{BEC-24},\cite{2022-6}.

Here, we describe discretizations confined to a bounded computational domain $\mathbb{D}$ under a homogeneous Dirichlet boundary condition. Spatial derivatives are approximated using the second-order finite difference (FD) method, and definite integrals are evaluated via the composite trapezoidal quadrature rule. For notational simplicity, we outline the detailed FD discretization procedure only for the one-dimensional case, as provided in \cite{2022-3,2022-35}. Extension to the two-dimensional case follows a similar approach. For the convergence analysis of this finite difference discretization toward the original energy functional optimization problem, we refer readers to \cite{2022-4}.

Let $\mathbb{D} = [a, b]$, and let $h$ be the spatial mesh size. Consider an equidistant partition of $\mathbb{D}$ given by $a = x_1 < x_2 < \ldots < x_{N-1} < x_N = b,$
with a total of $N$ partition points. Let $\phi_j$ denote the numerical approximation of $\phi(x_j)$ for $j = 1, 2, \ldots, N$, satisfying $\phi_1 = \phi(x_1) =0, \phi_N = \phi(x_N)= 0$
and define the vector of interior approximations as $\Phi = (\phi_2, \ldots, \phi_{N-1})^\top $. We have
 	$$
 	\begin{aligned}
 	E(\phi(x)) & \approx  \int_a^b \left[ \frac{1}{2} (\phi'(x))^2 + V(x)\phi(x)^2 + \frac{\beta}{2} \phi(x)^4 \right] \, {\rm d}x\\
 &=\sum_{j=1}^{N-1} \int_{x_j}^{x_{j+1}} \left[ -\frac{1}{2} \phi(x) \phi''(x) + V(x)\phi(x)^2 + \frac{\beta}{2} \phi(x)^4 \right] \, {\rm d}x\\& \approx
 	h \sum_{j=1}^{N-1} \left[ -\frac{1}{2} \phi_j \frac{\phi_{j+1} - 2\phi_j + \phi_{j-1}}{h^2} + V(x_j)\phi_j^2 + \frac{\beta}{2} \phi_j^4 \right]\\
 &=h \sum_{j=1}^{N-1}  \frac{1}{2} \left( \frac{\phi_{j+1} - \phi_j}{h} \right)^2 + h \sum_{j=2}^{N-1}V(x_j)\phi_j^2  + h \sum_{j=2}^{N-1}  \frac{\beta}{2} \phi_j^4 \\
 &=h \left[ \Phi^T B \Phi + \frac{\beta}{2} \sum_{j=2}^{N-1} \phi_j^4 \right],
 	\end{aligned}
 	$$
where $B = D + V \in \mathbb{R}^{(N - 2) \times (N - 2)} $  is a symmetric tri-diagonal matrix with entries
 	$$
 	D =
 	\begin{pmatrix}
 	\frac{1}{h^2} & - \frac{1}{2h^2}   &  &  &\\
 	- \frac{1}{2h^2}  & \frac{1}{h^2}  &\ddots  & & \\
 	& \ddots  & \ddots & \ddots &\\
 	&   & \ddots & \ddots &- \frac{1}{2h^2} \\
 	&   &  & - \frac{1}{2h^2} &\frac{1}{h^2}
 	\end{pmatrix},
 	\quad
 	V = \frac{1}{2}
 	\begin{pmatrix}
 	x_2^2 & \\
 	& \ddots &  \\
 	&  & x_{N-1}^2
 	\end{pmatrix}.
 	$$
Then the constraint with $d = 1$ can be truncated and discretized as
 	$$
 	\int_a^b \phi(x)^2 dx = \sum_{j=1}^{N-1} \int_{x_j}^{x_{j+1}} \phi(x)^2 {\rm d}x \approx h \sum_{j=2}^{N-1}  \phi_j^2 = 1.
    $$
Moreover, let $\mathbf{u}=\sqrt{h} \Phi$, we obtain the desired optimization problem
 	\begin{equation*}
 	\begin{aligned}
 	&\min_{\mathbf{u} \in \mathbb{R}^{N-2}} \hspace{0.5cm}\frac{\beta}{2h} \sum_{i=1}^{N-2} u_i^4 + \mathbf{u}^\top B \mathbf{u}\\
 	&  \hspace{0.4cm} \textup{s.t.}   \hspace{0.9cm}\|\mathbf{u}\|^2 = 1,
 	\end{aligned}
 	\end{equation*}
 	 where  $B$ is a symmetric positive definite sparse matrix.

Similarly, 	for the two-dimensional cases, the sizes of the variables are $(N - 2)^2$, and the matrix $B = D + V \in \mathbb{R}^{{(N - 2)}^2 \times {(N - 2)}^2} $  in the two-dimensional case has the following form:
 	 	\begin{equation*}
 	 	\begin{aligned}
 	 	&D_2 = \frac{1}{{h_x}^2}
 	 	\begin{pmatrix}
 	 	2 & -1  &  &  &\\
 	 	-1  & 2  & -1  & & \\
 	 	& \ddots  & \ddots & \ddots &\\
 	 	&   & -1 & 2 & -1 \\
 	 	&   &  & -1 & 2
 	 	\end{pmatrix},~h_x = \frac{b-a}{N-1},\\
 	 	& D = I_{N-2} \otimes D_2 + D_2 \otimes I_{N-2},~~~V = \frac{1}{2}
 	 	\begin{pmatrix}
 	 	{x_2}^2+{y_2}^2 & \\
 	 	& \ddots &  \\
 	 	&  & {x_{N-1}}^2  + {y_{N-1}}^2
 	 	\end{pmatrix}.
 	 	\end{aligned}
 	 	\end{equation*}
Unless otherwise specified, the computational domain $\mathbb{D}$ is set to $[-8, 8]$ in one-dimensional cases and
$[-8, 8] \times [-8, 8]$  in two-dimensional cases.

 \subsection{Alternating direction method of multipliers(ADMM)}
The problem (\ref{e1}) can be rewritten with the nonnegative constraint into the standard ADMM problem as follows:
\begin{equation*}
	\begin{aligned}
		&\min_{\mathbf{x} \in \mathbb{R}^n} \hspace{0.5cm} I_{\mathcal{S}}(\mathbf{x})+f(\mathbf{y})\\
		&  \hspace{0.2cm} \textup{s.t.}   \hspace{0.7cm}\mathbf{x}=\mathbf{y},
	\end{aligned}
\end{equation*}
where $I_{\mathcal{S}}(\mathbf{x})$ is the corresponding indicator function:
$$
I_{\mathcal{S}}(\mathbf{x}) = \begin{cases}
	0, &   \mathbf{x} \in \mathcal{S}, \\
	+\infty, & \mathbf{x} \notin \mathcal{S}.
\end{cases}
$$
and
$$
\mathcal{S} = \left\{  \mathbf{x} \,~|~ \| \mathbf{x} \| = 1,  \mathbf{x} \geq 0 \right\},~~f( \mathbf{y}) = \frac{\beta}{2h}\mathcal{T}_f \mathbf{y}^4+\mathbf{y}^\top B \mathbf{y}.
$$
For the reformulation above, assume the lagrange function is below:
$$
\mathcal{L}_{\rho}( \mathbf{x}, \mathbf{y}, \mathbf{w}) = I_{\mathcal{S}}( \mathbf{x})+f( \mathbf{y})+\mathbf{\mu}^\top( \mathbf{x}- \mathbf{y})+\frac{\rho}{2}\| \mathbf{x}- \mathbf{y} \, \|^2.
$$
Then we have the iteration steps for ADMM as follows:
$$\left\{
\begin{aligned}
	& \mathbf{x}^{k+1} := P_{\mathcal{S}}( \mathbf{y}^k-\frac{\mu^k}{\rho});\\
	& \mathbf{y}^{k+1} := {\arg\min} _{ \mathbf{y}} \left(f( \mathbf{y}) + \mu^{k^\top}( \mathbf{x}^{k+1}- \mathbf{y}) + \frac{\rho}{2}\| \mathbf{x}^{k+1}- \mathbf{y}\|^2 \right);\\
	&\mu^{k+1} := \mu^k + \rho( \mathbf{x}^{k+1}- \mathbf{y}^{k+1}),
\end{aligned}
\right.
$$
where $P_{\mathcal{S}}$ is the projection operator onto $\mathcal{S}.$ In the numerical experiments, we updates $ \mathbf{x}$ in the way $ \mathbf{x}^{k+1} := \frac{ \mathbf{y}^k-\frac{\mu^k}{\rho}}{\| \mathbf{y}^k-\frac{\mu^k}{\rho}\|}.$
As used in \cite{2022,2022-36}, we also invite the Riemannian gradient norm and set the stopping criterion for the ADMM as
$
\|\nabla f( \mathbf{x})- \mathbf{x}^\top \nabla f( \mathbf{x}) \mathbf{x}\| \leq 10^{-6}.
$
The solver for the convex subproblem in ADMM is the Newton method, and furthermore, the Gauss-Seidel method is used to get a descent direction for the Newton method. Unless otherwise specified, the stopping criterion in the $k$-th outer iteration for the Newton method solving subproblems is $\frac{newDec^2}{2} \leq \epsilon_k$, where $newDec$ is the Newton decrement and $\epsilon_k = 10^{-2} \cdot 0.5^k.$
We observed that if $\rho$ is chosen too small, the algorithm might be divergent. In general, the larger the interaction coefficient $\beta$, the larger an appropriate $\rho$ will be needed.

 \subsection{Implementation details}
We consider the fourth degree polynomial optimization problem (\ref{e3}) where the underlying tensor $\mathcal{T}_f$ is generated by (\ref{e2}). To investigate the numerical behaviors of our proposed Algorithm  \ref{alg1} (denoted by PAM), we compare it with the state-of-the-art algorithm, i.e., ADMM \cite{2022,2022--1,2022-32}.
In our algorithmic implementation, we take
\begin{equation*}
\mathrm{Err}:= \frac{| f_{\alpha}(\mathbf{x}^{(k+1)})- f_{\alpha}(\mathbf{x}^{(k)}) |}{\max \{ | f_{\alpha}(\mathbf{x}^{(k+1)}) |,| f_{\alpha}(\mathbf{x}^{(k)}) |, 1 \} }  \leq 10^{-6}
\end{equation*}
as the stopping criterion for our method and set the maximum iteration as 2000 for all methods. We report the number of iterations (denoted by $\mathbf{Total~iter}$), computing time in seconds (denoted by $\mathbf{cpu}$), and the final objective function values (denoted by $\mathbf{obj~val}$).
All numerical experiments were carried out in MATLAB R2018b on a desktop computer equipped with an Intel(R) Core(TM) i7-10700 CPU @ 2.90GHz, 16 GB of RAM, running Microsoft Windows 10 (64-bit).

\subsection{Comparison between the PAM and ADMM method}
\begin{figure}[htbp]
	\centering
	\captionsetup[subfigure]{skip=2pt}
\begin{subfigure}{0.48\textwidth}
	\centering
	\includegraphics[width=5.8cm, trim=30 200 50 210, clip]{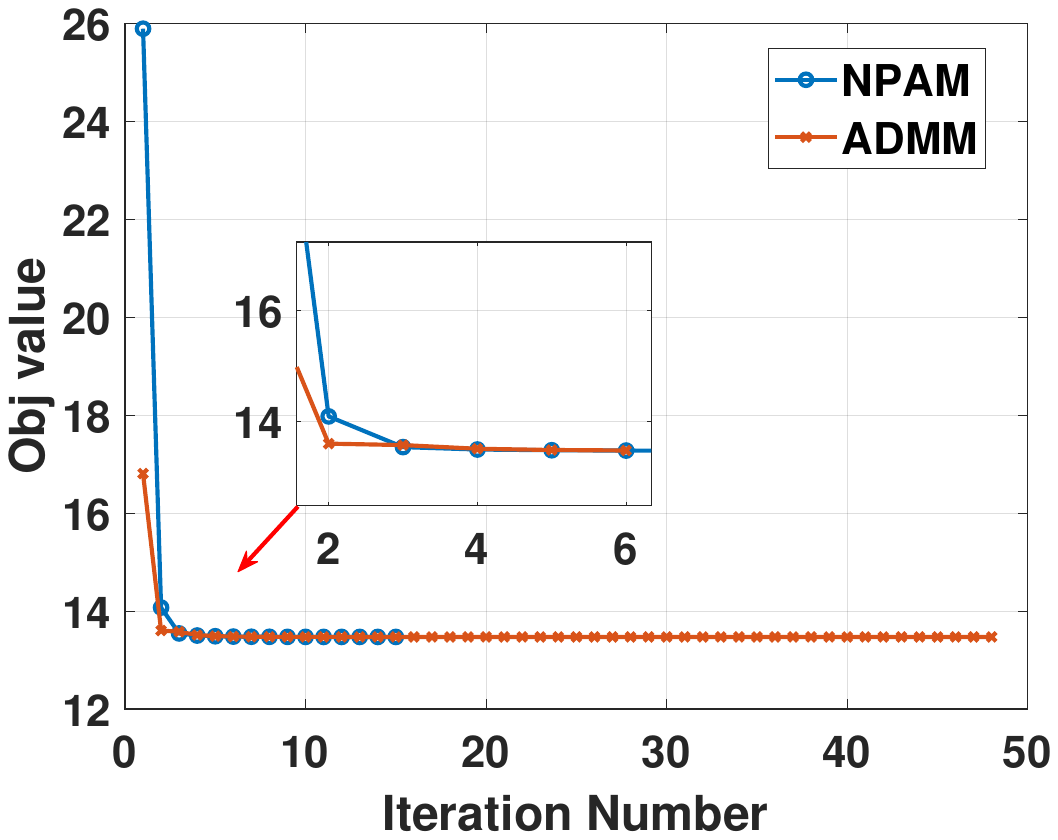}
	\caption{$\beta = 250$}
\end{subfigure}
	\hfill
	\begin{subfigure}{0.48\textwidth}
		\centering
	\includegraphics[width=6cm, trim=30 200 50 210, clip]{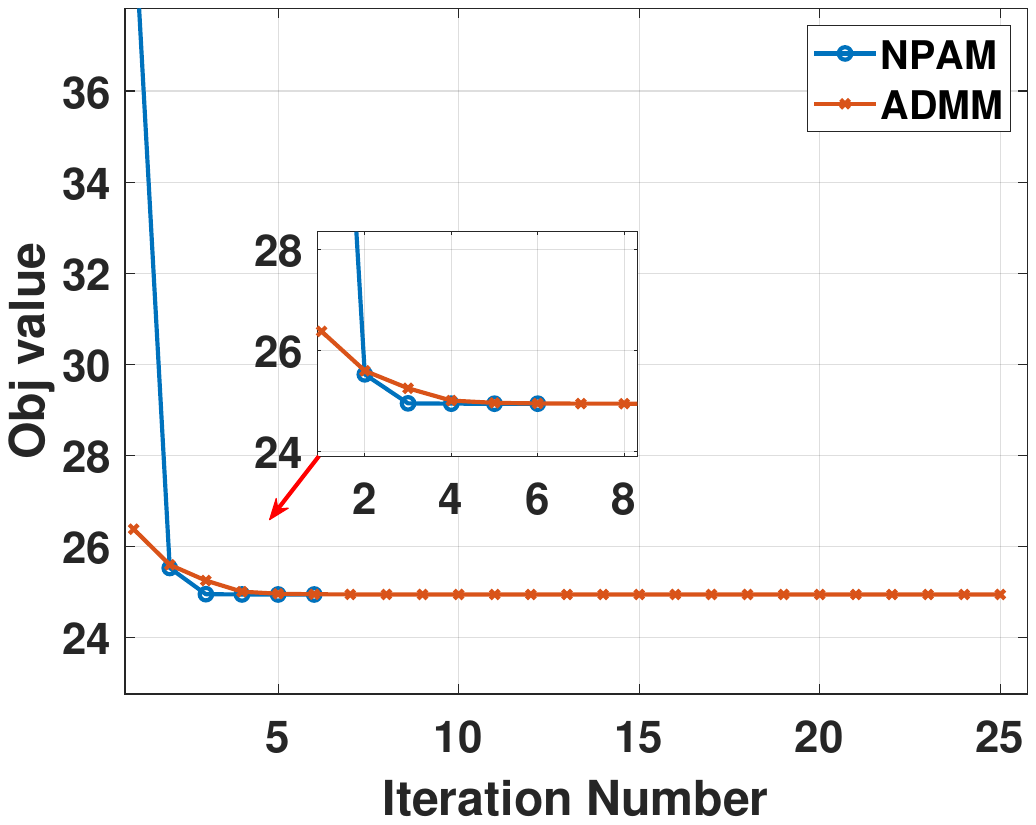}
		\vspace{-3pt}
		\caption{$\beta = 500$}
	\end{subfigure}
	\vskip\baselineskip
	\begin{subfigure}{0.48\textwidth}
		\centering
		\includegraphics[width=6cm, trim=30 200 50 210, clip]{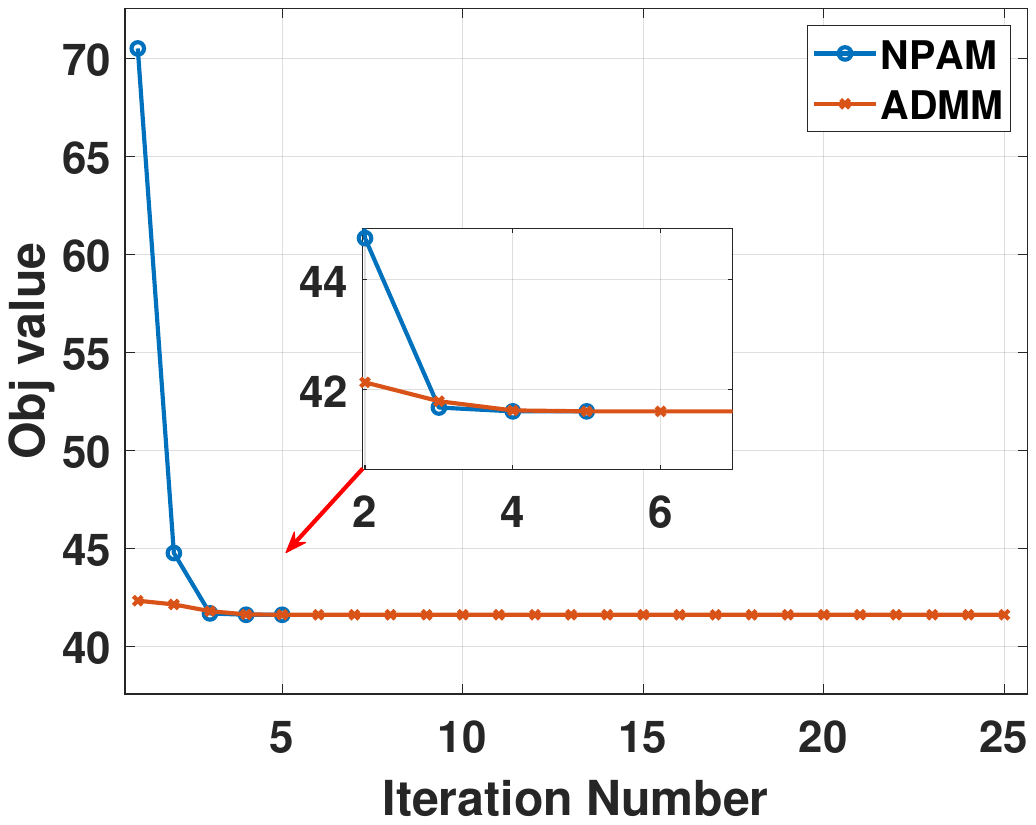}
		\vspace{-3pt}
		\caption{$\beta = 1000$}
	\end{subfigure}
	\hfill
	\caption{Illustration of the objective function value versus total iteration numbers for PAM and ADMM in the one-dimensional case with $N=50$ and different values of $\beta$.}
	\label{fig:1}
\end{figure}

\begin{figure}[H]
	\centering
	\captionsetup[subfigure]{skip=2pt}
	\begin{subfigure}{0.48\textwidth}
		\centering
		\includegraphics[width=6cm, trim=30 200 50 210, clip]{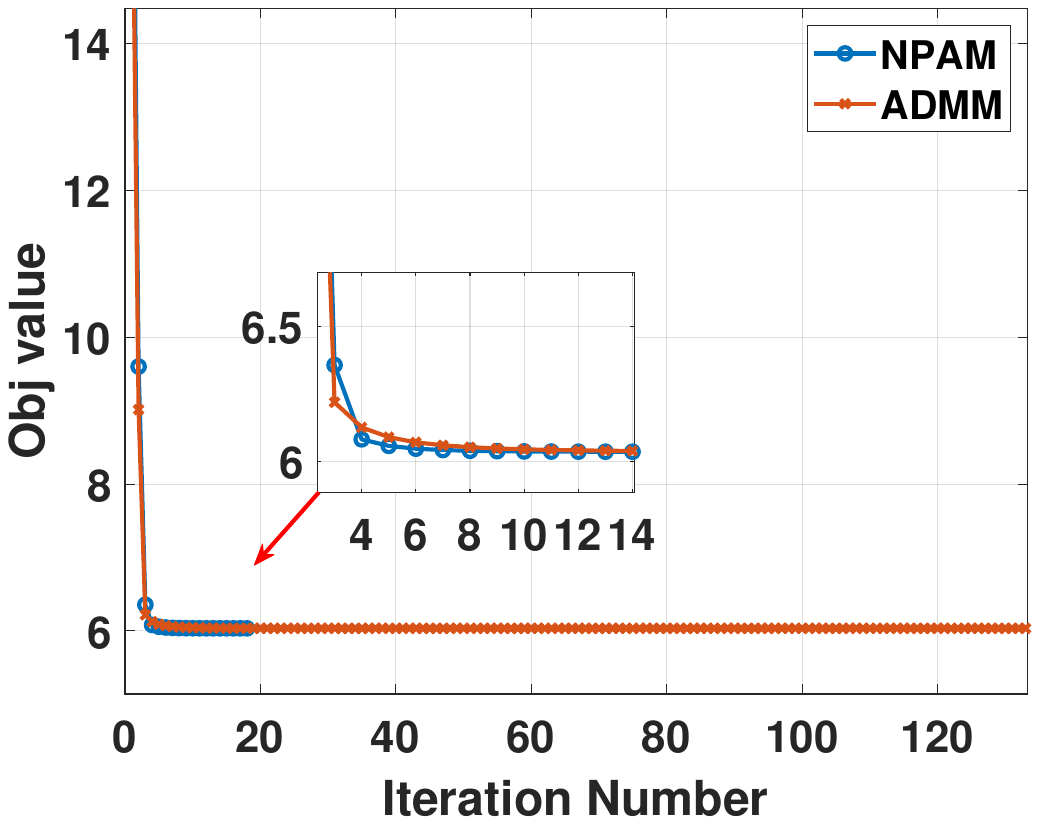}
		\caption{$\beta = 250$}
	\end{subfigure}
	\hfill
	\begin{subfigure}{0.48\textwidth}
		\centering
		\includegraphics[width=6cm, trim=30 200 50 210, clip]{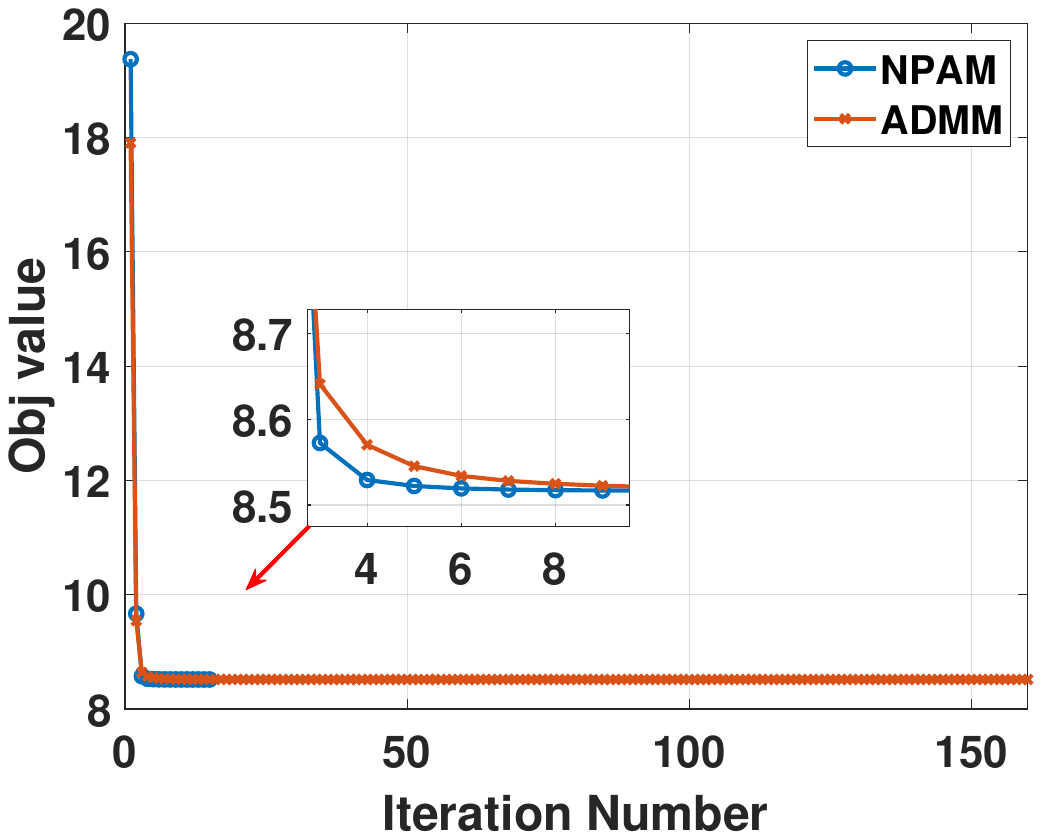}
		\vspace{-3pt}
		\caption{$\beta = 500$}
	\end{subfigure}
	\vskip\baselineskip
	\begin{subfigure}{0.48\textwidth}
		\centering
		\includegraphics[width=6cm, trim=5 180 15 210, clip]{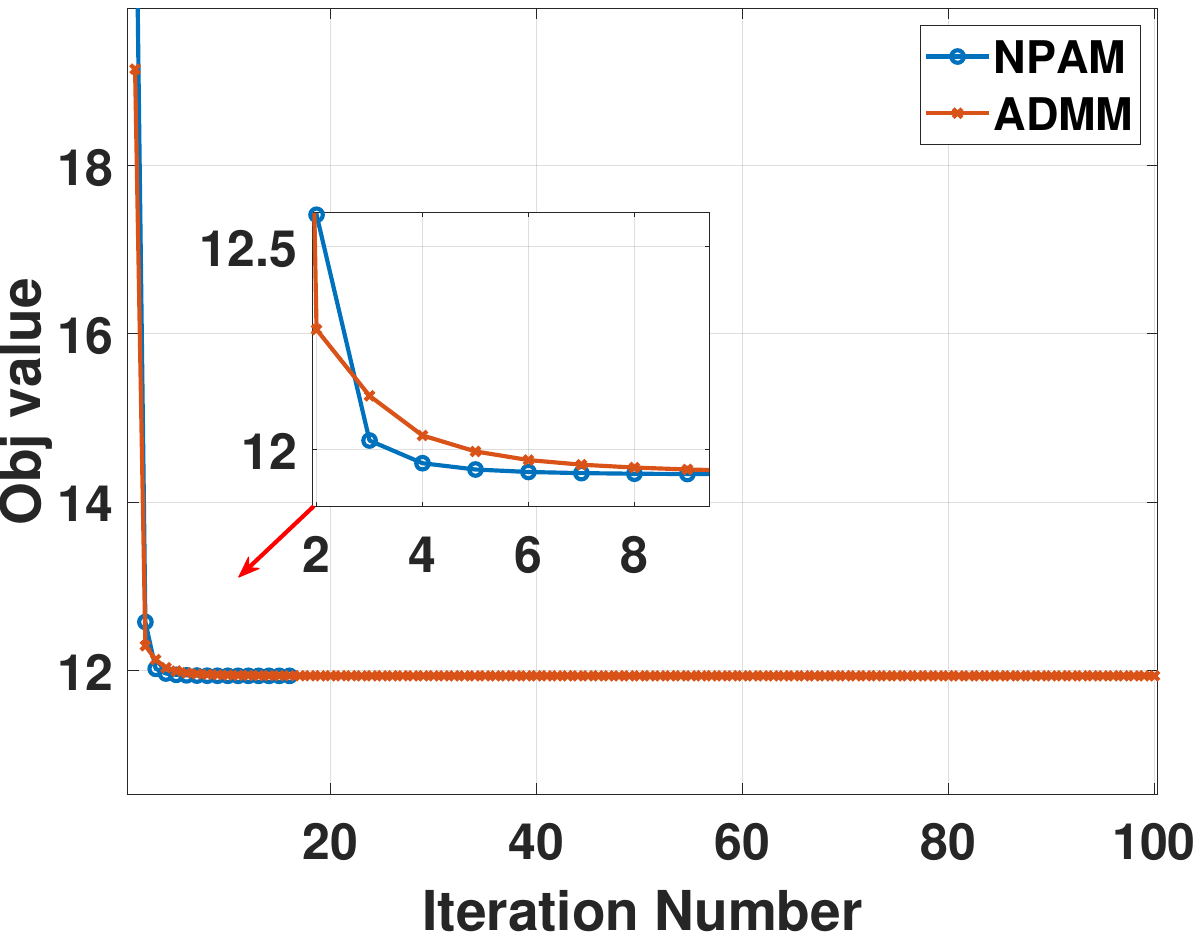}
		\vspace{-1pt}
		\caption{$\beta = 1000$}
	\end{subfigure}
	\hfill
	\caption{Illustration of the objective function value versus total iteration numbers for PAM and ADMM in the two-dimensional case with $N=10$ and different values of $\beta$.}
	\label{fig:2}
\end{figure}

We compare PAM with the ADMM method for quartic--quadratic optimization problems in both one-dimensional and two-dimensional settings. In the one-dimensional case, we test $\beta = 250, 500, 1000$ with $N = 10,20,30,40,50, 60,70$, while in the two-dimensional case we consider the same $\beta$ values with grid sizes $N = 7,8,9,10,11$. The parameters of PAM are fixed as $\gamma_i = 0.5$ for $i = 1,2,3,4$, and all methods are initialized with the same normally distributed random vectors. The numerical results in Tables \ref{table vs 1d} and \ref{table vs 2d} demonstrate that PAM consistently requires fewer iterations and less CPU time than ADMM, while achieving comparable objective values. To further illustrate convergence behavior, Fig.~\ref{fig:1} (1D, $N=50$) and Fig.~\ref{fig:2} (2D, $N=10$) present representative cases with $\beta = 250,500,1000$. In each figure, the left panels depict the overall convergence profiles, while the right panels highlight the rapid initial convergence of PAM.

As shown in Table \ref{table vs 1d}, PAM achieves nearly identical objective values to ADMM but with substantially fewer iterations and reduced computational time. In most cases, ADMM exhibits a faster decrease in the objective function during the initial iterations. However, as the iteration proceeds, PAM ultimately requires fewer total iterations to reach a solution of comparable accuracy. A possible explanation is that ADMM involves inner Newton iterations, where the descent direction is computed only approximately by solving a linear system to very low precision. This low-accuracy solution may yield a rapid decrease at the beginning but limits the efficiency of ADMM as it approaches the optimal value. This phenomenon may warrant further investigation.

\begin{table}[H]
	\fontsize{6.5}{9}\selectfont
		\caption{Comparison between PAM and ADMM method in the one-dimensional case. The columns of total iter show the number of iterations. For ADMM method, it includes the iteration of Newton method, and the numbers in brackets stand for the outer ADMM iteration.}
	\label{table vs 1d}
	\centering
	\begin{threeparttable}
		\resizebox{0.9\textwidth}{!}{     
		\begin{tabular}{cccccccccc}
        	\toprule               
			\multirow{2}{*}{\textbf{$\beta$}}          
           &\multirow{2}{*}{$N$}                   
           &&\multicolumn{3}{c}{\textbf{PAM}}	         
           	&&\multicolumn{3}{c}{\textbf{ADMM}}\\	
           	\cmidrule(l){4-6}                                                   
           	\cmidrule(l){8-10}
           	&&&Total iter&cpu(s)&obj val&&Total iter&cpu(s)&obj val\\
           	
			\midrule      
		\multirow{4}{*}{250($\alpha$=15, $\rho=80$)}
		&10&&5.6&0.0006&15.6145&&50(23)&0.0106&15.6145\\
		&20&&11.8&0.0025&15.6246&&72(34)&0.0429&15.6246\\
		&30&&11.8&0.0039&15.6247&&62(29)&0.0286&15.6246\\
		&40&&11.2&0.0076&15.6248&&59(28)&0.0302&15.6248\\
		&50&&11.1&0.0196&15.6249&&65(31)&0.0406&15.6249\\
		&60&&11.1&0.0259&15.6249&&65(31)&0.0372&15.6249\\
		&70&&10.8&0.0421&15.6249&&67(32)&0.0792&15.6249\\				
		\midrule      
		\multirow{4}{*}{500($\alpha$=27, $\rho=115$)}
		&10&&4.1&0.0006&25.1569&&78(37)&0.0168&25.1569\\
		&20&&4.6&0.0013&24.9847&&54(25)&0.0279&24.9847\\
		&30&&5.1&0.0024&24.9536&&54(25)&0.0212&24.9536\\
		&40&&5.4&0.0034&24.9449&&54(25)&0.0232&24.9450\\
		&50&&5.4&0.0103&24.9428&&54(25)&0.0322&24.9428\\
		&60&&5.3&0.0188&24.9427&&54(25)&0.0694&24.9427\\
		&70&&5.8&0.0284&24.9431&&54(25)&0.0890&24.9431\\
		\midrule      
		\multirow{4}{*}{1000($\alpha$=50,, $\rho=210$)}
		&10&&4.1&0.0006&43.1132&&62(29)&0.0103&43.1132\\
		&20&&4.6&0.0008&42.0423&&56(26)&0.0114&42.0423\\
		&30&&5.2&0.0022&41.7744&&56(26)&0.0168&41.7744\\
		&40&&4.5&0.0036&41.6631&&54(25)&0.0123&41.6631\\
		&50&&4.3&0.0092&41.6080&&54(25)&0.0394&41.6080\\
		&60&&4.5&0.0144&41.5787&&54(25)&0.0458&41.5787\\
		&70&&4.6&0.0242&41.5628&&54(25)&0.0693&41.5628\\
		\bottomrule     
		\end{tabular}
	}
	\end{threeparttable}
\end{table}

\begin{table}[H]
	\fontsize{6.5}{9}\selectfont
	\caption{Comparison between PAM and ADMM method in the two-dimensional case. The columns of total iter show the number of iterations. For ADMM method, it includes the iteration of Newton method, and the numbers in brackets stand for the outer ADMM iteration.}
	\label{table vs 2d}
	\centering
	\begin{threeparttable}
		\resizebox{0.9\textwidth}{!}{
			
			\begin{tabular}{cccccccccc}
				\toprule               
				\multirow{2}{*}{\textbf{$\beta$}}          
				&\multirow{2}{*}{$N$}                   
				&&\multicolumn{3}{c}{\textbf{PAM}}	         
				&&\multicolumn{3}{c}{\textbf{ADMM}}\\	
				\cmidrule(l){4-6}                                                   
				\cmidrule(l){8-10}
				&&&Total iter&cpu(s)&obj val&&Total iter&cpu(s)&obj val\\
				
				\midrule      
				\multirow{4}{*}{250($\alpha$=6,$\rho$=40)}
				&7&&10.3&0.0042&6.0594&&127(81)&0.0938&6.0594\\
				&8&&12.1&0.0082&6.1161&&240(167)&0.1215&6.1161\\
				&9&&18.1&0.0288&6.0696&&320(194)&0.1188&6.0696\\
				&10&&16.9&0.0490&6.0355&&316(198)&0.1815&6.0355\\
				&11&&13.3&0.0863&6.0533&&257(161)&0.1812&6.0533\\
				
				\midrule      
				\multirow{4}{*}{500($\alpha$=7,$\rho$=50)}
				&7&&16.8&0.0038&8.3697&&311(183)&0.1615&8.3697\\
				&8&&10.8&0.0055&8.4589&&193(109)&0.0978&8.4589\\
				&9&&10.1&0.0147&8.5092&&181(111)&0.0706&8.5092\\
				&10&&14.8&0.0360&8.5164&&281(160)&0.1226&8.5085\\
				&11&&13.1&0.0717&8.5085&&266(155)&0.1315&8.5085\\
				
				\midrule      
				\multirow{4}{*}{1000($\alpha$=9,$\rho$=65)}
				&7&&13.9&0.0044&12.0468&&234(120)&0.0715&12.0468\\
				&8&&11.8&0.0061&11.9693&&163(84)&0.0725&11.9693\\
				&9&&16.2&0.0244&11.9471&&362(193)&0.1058&11.9471\\
				&10&&14.3&0.0429&11.9391&&245(128)&0.0841&11.9391\\
				&11&&11.2&0.0739&11.9566&&178(91)&0.0845&11.9566\\
				\bottomrule     
			\end{tabular}
		}
	\end{threeparttable}
\end{table}

\subsection{Ground state profiles}
\begin{figure}[H]
 	   \vspace{-6em} 
	\begin{minipage}[t]{0.46\linewidth} 
		\centering
		\includegraphics[width=8cm]{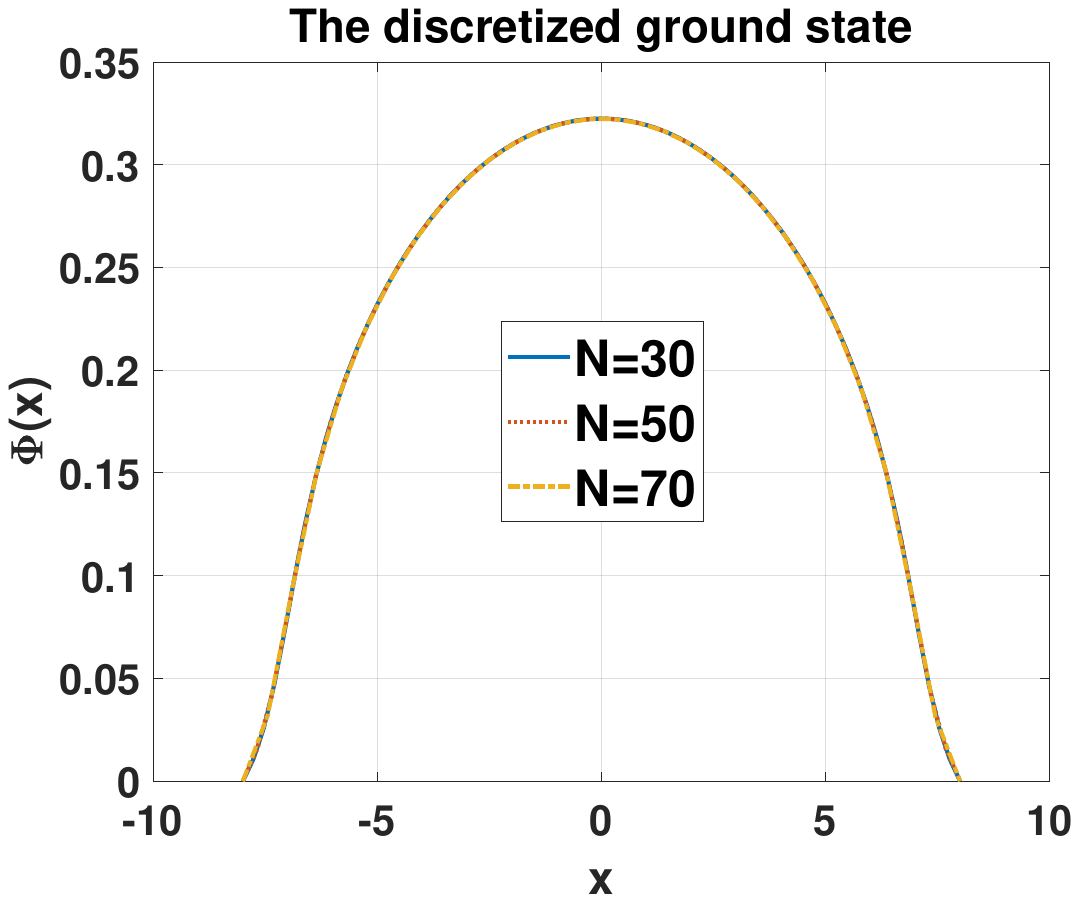}
	\end{minipage}
	\begin{minipage}[t]{0.5\linewidth} 
		\hspace{2mm}
		\includegraphics[width=8.5cm]{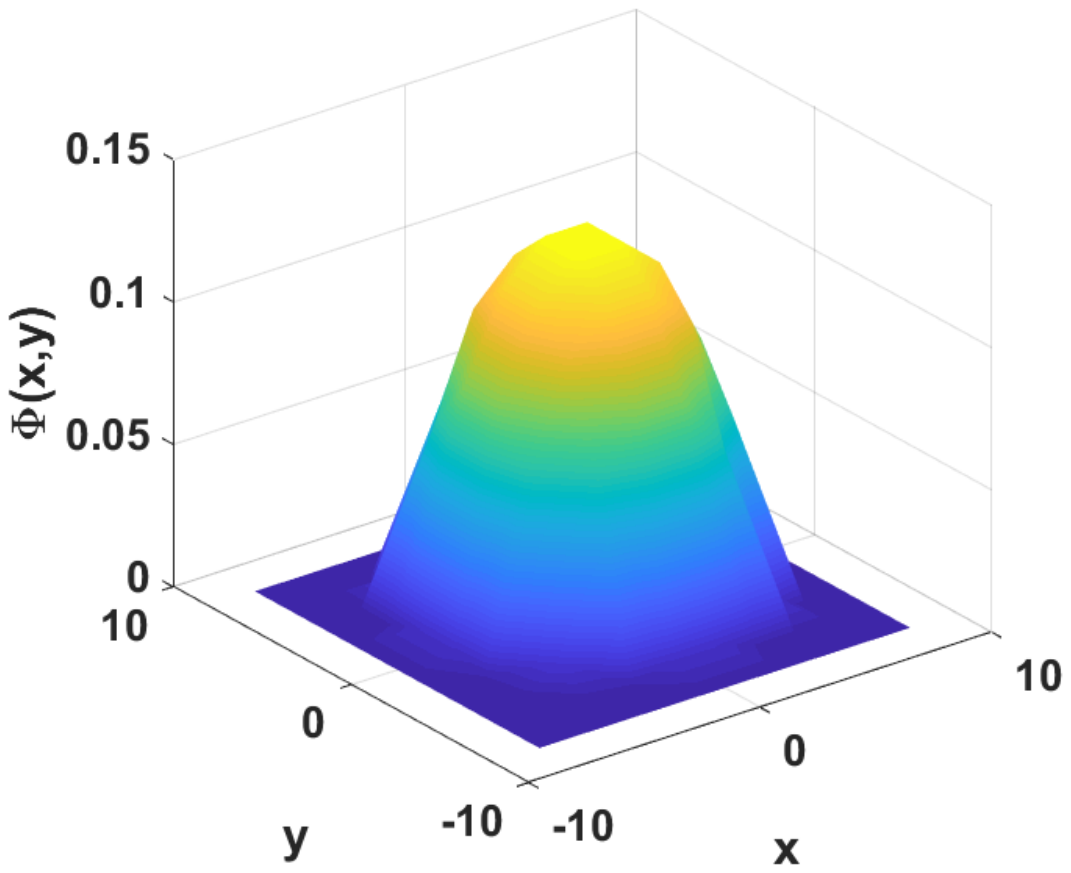}
	\end{minipage}
  \vspace{-7em} 
	\caption{The discretized ground state computed by PAM. The left one is for the one-dimensional space with $N$ = 30, 50, 70. The right one is for the two-dimensional space with $N$ = 10.}\label{fig7}
\end{figure}

Fig.~\ref{fig7} shows some examples of the discretized ground state computed by PAM method in the one-dimensional and two-dimensional cases. The obtained solutions are nonnegative, smooth, and symmetric, with rapid decay at the boundary, which coincides with the expected physical properties of the BEC ground state. These results indicate that PAM can stably capture the physically meaningful ground state in different spatial dimensions.

\subsection{Sensitivity analysis of $\alpha$}
\begin{figure}[H]
   	    \vspace{-0.5em} 
	\begin{minipage}[t]{0.5\linewidth} 
		\centering
		\includegraphics[width=7cm]{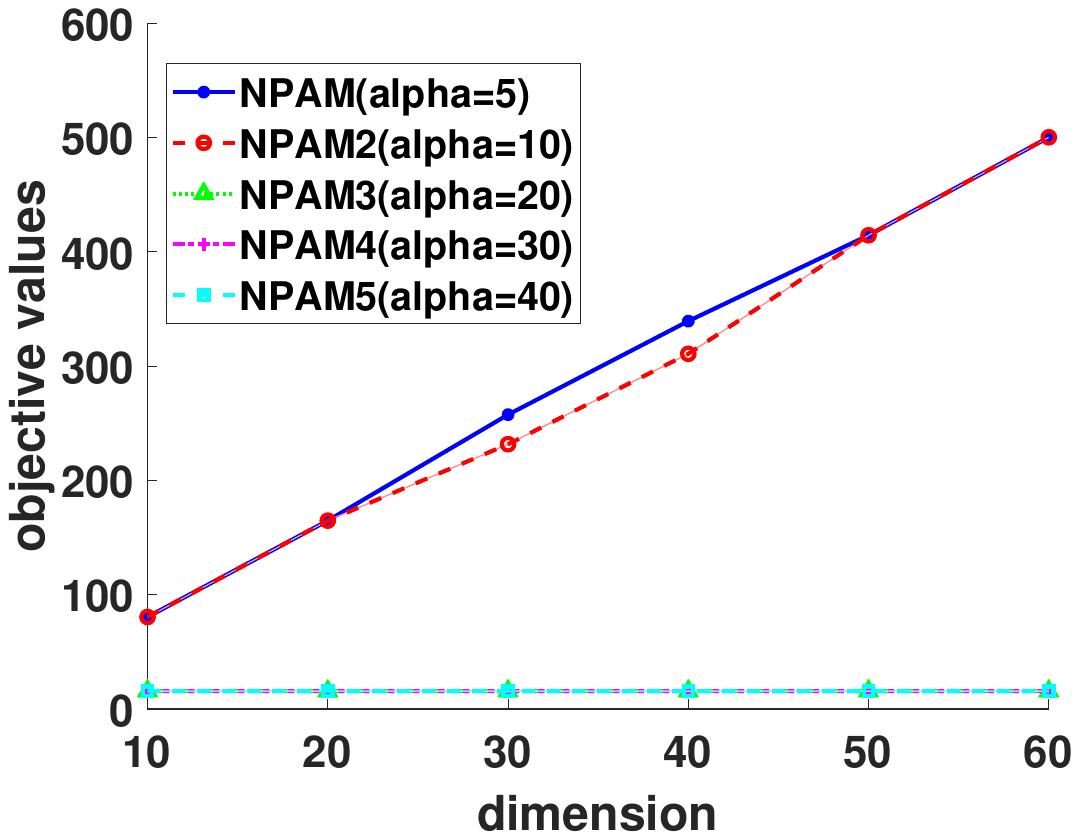}
	\end{minipage}
	\begin{minipage}[t]{0.5\linewidth} 
		\hspace{2mm}
		\includegraphics[width=7cm]{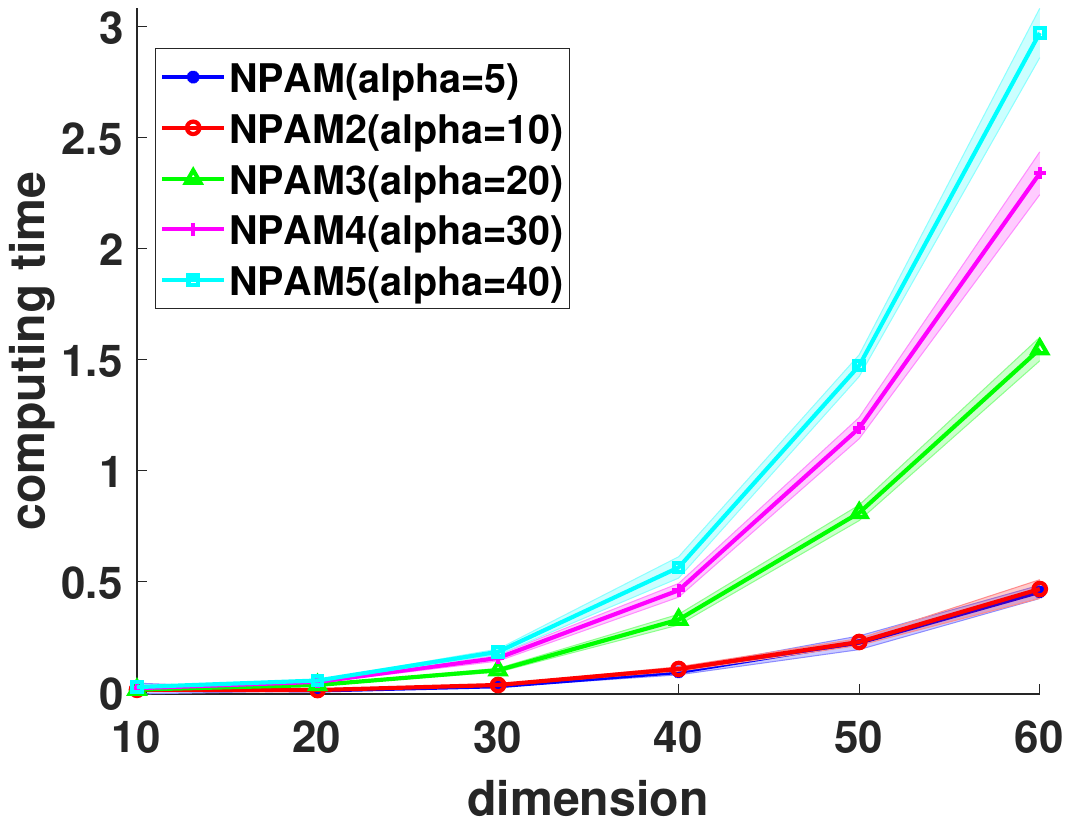}
	\end{minipage}
  \vspace{-5em} 
	\caption{Comparison between different $\alpha$ with PAM method in the one-dimensional case for $ N = 10,20,30,40,50,60 $and $\beta=250$. The left one show the objective values, and the right one stand for the computing time.}\label{fig8}
	 	\vspace{-1em} 
\end{figure}

\begin{figure}[H]
	    \vspace{-0.5em} 
  	\begin{minipage}[t]{0.5\linewidth} 
  		\centering
  		\includegraphics[width=7cm]{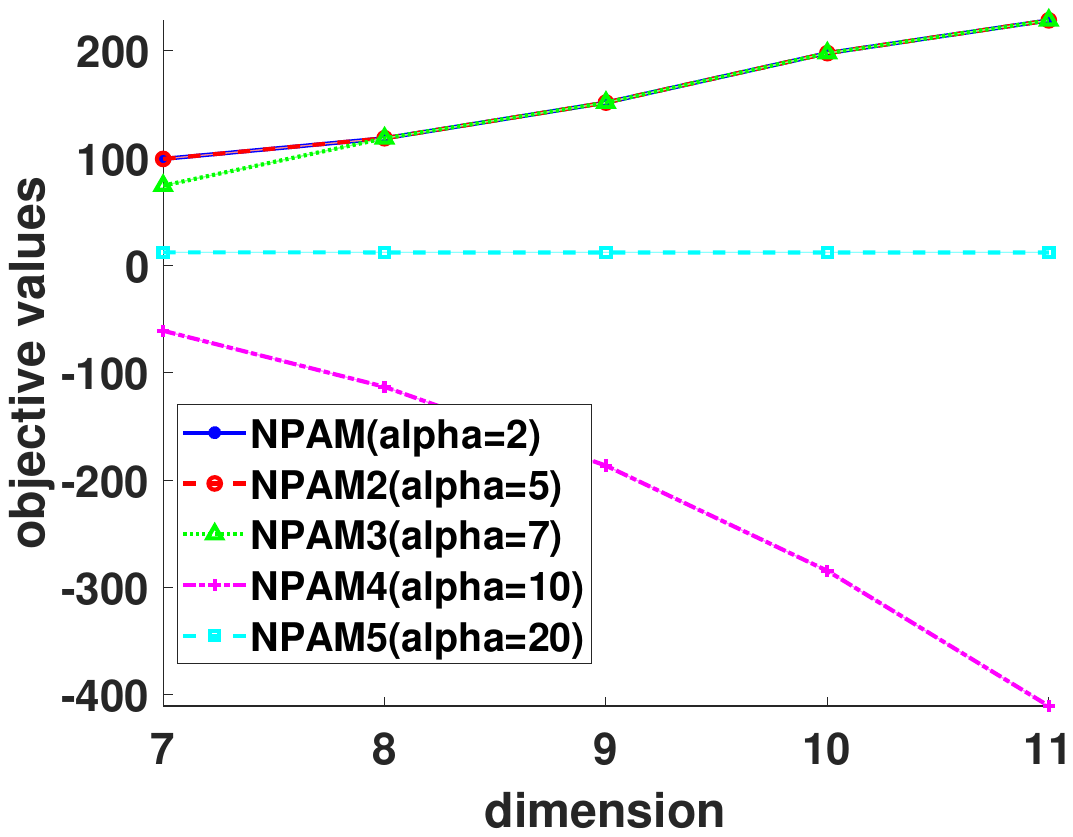}
  	\end{minipage}
  	\begin{minipage}[t]{0.5\linewidth} 
  		\hspace{2mm}
  		\includegraphics[width=7cm]{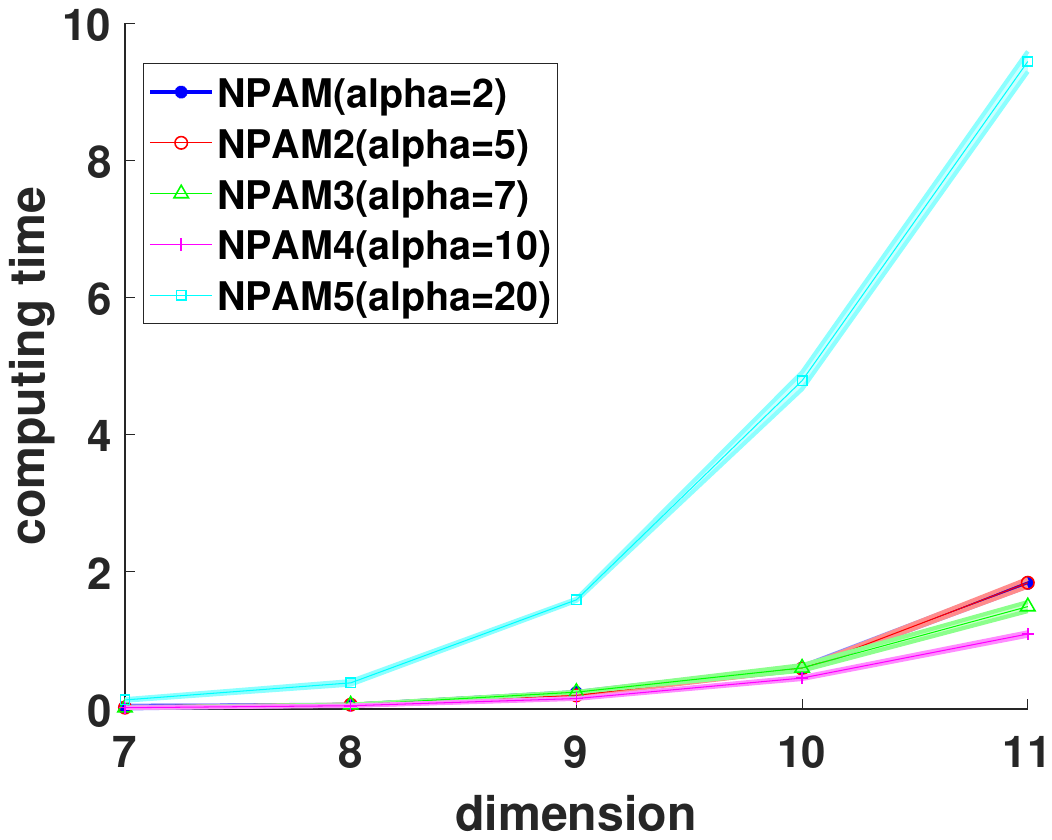}
  	\end{minipage}
  \vspace{-5em} 
  	\caption{Comparison between different $\alpha$ with PAM method in the two-dimensional case for $N = 7,8,9,10,11$ and $\beta=500$. The left one show the objective values, and the right one stand for the computing time.}\label{fig9}
  \end{figure}

Note that augmented (or regularized) terms are introduced in equations~(\ref{e9}) and~(\ref{e10}). We next examine the sensitivity of the parameter $\alpha$. In the one-dimensional case (Fig.~\ref{fig8}), we set $\beta = 250$, $N = 10,20,30,40,50,60$, and test $\alpha=5,10,20,30,40$. In the two-dimensional case (Fig.~\ref{fig9}), we set $\beta = 500$, $N =7,8,9,10,11 $ and test $\alpha = 2,5,7,10,20 $.
In both figures, the left panels report the averaged objective values, while the right panels display the computational times. Solid lines indicate the mean over five runs, and shaded areas represent the standard deviation. Since $N$ denotes the number of grid points, the optimal objective values are expected to remain comparable across different $N$. Nevertheless, the upward trend in Fig.~\ref{fig8} suggests instability for certain $\alpha$, while the right panels show that appropriately chosen $\alpha$ can reduce running time under stable performance. Similar behavior is observed in the two-dimensional case (Fig.~\ref{fig9}).

Overall, the method exhibits robust performance across a broad range of $\alpha$. Extremely small $\alpha$ may impair stability, whereas excessively large $\alpha$ increase computational time. In practice, a moderate choice of $\alpha$ achieves a favorable balance between stability and efficiency.

\medskip
\begin{acknowledgements}
This work was supported by National Natural Science Foundation of China (12071249), Shandong Provincial Natural Science Foundation for Distinguished Young Scholars (ZR2021JQ01), Shandong Provincial Natural Science Foundation (ZR2024MA003).
\end{acknowledgements}

\medskip
\noindent{\bf Declarations}
The authors declare that they have no conflict of interest.

\medskip
\noindent{\bf Data Availability Statement}
The data that support the findings of this study is available from the corresponding author upon reasonable request.


\end{document}